\documentclass[brochure,12pt]{bourbaki}
\usepackage[matrix,arrow]{xy}
\usepackage{amssymb,footnote,sabbah_Exp1050}
\usepackage[francais]{babel}
\date{Janvier 2012}
\bbkannee{64\`eme ann\'ee, 2011-2012}
\bbknumero{1050}
\title[Théorie de Hodge sauvage]{Théorie de Hodge et correspondance de~Hitchin-Kobayashi sauvages}
\subtitle{d'apr\`es T. Mochizuki}
\author{Claude SABBAH}
\address{Centre de Math\'ematiques Laurent Schwartz\\
\'Ecole Polytechnique\\
91128 Palaiseau Cedex}
\email{sabbah@math.polytechnique.fr}
\urladdr{http://www.math.polytechnique.fr/~sabbah}
\thanks{Cette recherche a été effectuée dans le cadre du programme  ANR-08-BLAN-0317-01 de l'Agence nationale de la recherche.}

\begin{document}
\frontmatter
\maketitle
\tableofcontents
\mainmatter

\section*{Introduction: le théorème de Lefschetz difficile}
\subsubsection*{Les $\cD$-modules holonomes simples}
Soient~$Z^o$ une variété quasi-projective lisse complexe irréductible et $(V,\nabla)$ un fibré vectoriel algébrique muni d'une connexion intégrable (\ie sans courbure), sans sous-fibré propre stable par la connexion (\ie simple). C'est le type d'objet auquel on s'intéresse dans cet exposé. Choisissons une compactification projective~$j:Z^o\hto Z$ et un plongement de~$Z$ dans une variété projective lisse~$X$. Il est connu qu'un tel $(V,\nabla)$ se prolonge de manière unique en un $\cD_X$-module (\ie $\cO_X$-module avec connexion intégrable $\nabla$) holonome simple à support dans~$Z$ et qu'on obtient ainsi tous les $\cD_X$-modules holonomes simples à support dans~$Z$, lisses sur $Z^o$.

Pour un tel $\cD_X$-module $\cM$, le complexe de de~Rham analytique
\[
\DR\cM:=(\Omega_{X^\an}^{\cbbullet+\dim X}\otimes\cM,\nabla)
\]
est à cohomologie $\CC$-constructible, d'après un théorème de Kashiwara. Plus précisément, c'est un \emph{faisceau pervers}.

Quels sont les faisceaux pervers qu'on obtient de cette manière? On ne connaît pas la réponse à cette question. Néanmoins, on sait \cite{Kashiwara84,Mebkhout84,Mebkhout84b} que la correspondance de Riemann-Hilbert $\cM\mto\DR\cM$ est une équivalence entre la catégorie des $\cD_X$-modules holonomes réguliers (analytiques) et celle des faisceaux pervers. Ceci, joint au théorème GAGA pour les modules holonomes réguliers \cite{K-K81}, implique que tout \emph{faisceau pervers simple} est obtenu de cette manière. D'après \cite{B-B-D81} et \cite{G-M83bis}, un tel faisceau pervers n'est autre (au décalage par $\dim Z$ près) que le complexe d'intersection $\IC_Z(\cL)$ de Goresky-MacPherson (extension intermédiaire $j_{!*}\cL$) associé à un faisceau localement constant irréductible~$\cL$ sur~$Z^o$ (\ie une représentation linéaire irréductible de $\pi_1(Z^o,\star)$).

Il est par ailleurs facile d'obtenir des complexes $\DR\cM$ qui ne sont pas des faisceaux pervers simples, ni même semi-simples (\ie sommes directes d'objets simples), si on accepte que $\cM$ simple ait des singularités irrégulières. Une façon d'obtenir de tels exemples consiste à utiliser la transformation de Fourier.

\begin{exem*}
Soient $T_1,\dots, T_r$ ($r\geq2$) des éléments de $\GL_n(\CC)$ ($n\geq2$) dont le produit est égal à l'identité et qui n'ont pas de vecteur propre commun. Supposons aussi que~$1$ ne soit pas valeur propre des $T_i$ (multiplier chaque $T_i$ par $\lambda_i\in\CC^*$ assez général, en faisant en sorte que $\prod\lambda_i=1$). Ils définissent une représentation irréductible du groupe fondamental de $\PP^1$ privé de $r$ points, tous à distance finie, donc un faisceau localement constant irréductible de rang~$n$ sur cet espace. Son complexe d'intersection sur $\Afu$ (coordonnée $z$) correspond à un module holonome~$M$ sur l'algèbre de Weyl $\CC[z]\langle\partial_z\rangle$, dont toutes les singularités sont régulières.

Le transformé de Fourier $\Fou M$ est $M$ lui-même sur lequel on voit $\partial_z$ opérer comme la multiplication par une variable $\zeta$ et $-z$ comme la dérivation $\partial_\zeta$. Si $M$ est simple, $\Fou M$ l'est aussi, mais on peut montrer (voir par exemple \cite{Malgrange91}) que ce dernier a une singularité irrégulière en $\zeta=\infty$, une singularité régulière en $\zeta=0$, et pas d'autre singularité. On peut aussi montrer que les hypothèses faites sur les~$T_i$ impliquent que, sur l'ouvert $\zeta\neq0$, le faisceau localement constant de ses sections horizontales est de rang égal à $nr$, et sa monodromie a pour seule valeur propre $1$, avec $r$ blocs de Jordan de taille $2$ et $(n-2)r$ blocs de taille $1$. Ce faisceau localement constant n'est donc pas semi-simple.

Soit $\FcM$ l'unique $\cD_{\PP^1}$-module holonome simple dont la restriction à~$\Afu$ (coordonnée~$\zeta$) est~$\Fou M$. Si $\DR\FcM$ était pervers semi-simple, il serait somme directe de faisceaux à support ponctuel et de complexes d'intersection de faisceaux localement constants irréductibles (à~un décalage près), parmi lesquels le faisceau localement constant ci-dessus, d'où une contradiction.
\end{exem*}

Nous allons voir cependant que les faisceaux pervers $\DR\cM$, pour $\cM$ simple, satisfont tous au théorème de Lefschetz difficile. Dans la suite, nous travaillerons dans le cadre analytique complexe uniquement, contrairement au début de cette introduction.
 
\subsubsection*{Le théorème de Lefschetz difficile}
Dans différents exposés en 1996 (voir \cite{Kashiwara98b}), Kashiwara  a conjecturé une version très générale du théorème de Lefschetz difficile en géométrie algébrique complexe, qui a été démontrée récemment par T\ptbl Mochizuki \cite{Mochizuki08}:

\begin{theo}\label{th:HLT}
Soit~$X$ une variété algébrique projective complexe lisse et soit $L$ l'opérateur de cup-produit par la classe de Chern d'un fibré en droites ample sur~$X$. Alors, pour tout $\cD_X$-module holonome simple $\cM$ sur~$X$ et tout $k\geq1$, l'itéré $k$-ième $L^k:\bH^{-k}(X,\DR\cM)\to\bH^k(X,\DR\cM)$ est un isomorphisme.
\end{theo}

\begin{rema}\label{rem:H0}
Avec le décalage définissant $\DR$, $\bH^k$ peut être non nul seulement pour $k\in[-\dim X,\dim X]$. De plus, pour $\cM$ simple et différent de la connexion triviale $(\cO_X,\rd)$, on a aussi nullité pour $k=-\dim X,\dim X$. En effet, considérons d'abord le complexe $\Gamma(X,\DR\cM)$: son $H^{-\dim X}$ est l'espace des sections de $\cM$ sur~$X$ annulées par~$\nabla$; il y a donc un sous-module $H^{-\dim X}\otimes(\cO_X,\rd)$ contenu dans $\cM$, et l'hypothèse \og $\cM$ simple et $\neq(\cO_X,\rd)$\fg implique $H^{-\dim X}=0$. La suite spectrale d'hypercohomologie montre que $\bH^{-\dim X}\subset H^{-\dim X}$, d'où la nullité de $\bH^{-\dim X}$. Un argument de dualité permet aussi d'en déduire $\bH^{\dim X}=0$. Ainsi, l'énoncé \ref{th:HLT} est de peu d'intérêt lorsque~$X$ est une courbe.
\end{rema}

Pour obtenir un théorème de ce type, on montre d'abord l'existence d'une structure plus riche, qui fait apparaître une notion de pureté \cite{DeligneHI}. La théorie de Hodge joue ce rôle en géométrie algébrique complexe \og modérée\fg (voir l'excellent panorama \cite{C-M09}):
\begin{enumeratea}
\item\label{enum:caslisse}
Si $(V,\nabla)$ est un fibré holomorphe à connexion intégrable sur~$X$, le complexe $\DR(V,\nabla)=(\Omega_X^{\cbbullet+\dim X}\otimes V,\nabla)$ n'a de cohomologie qu'en degré $-\dim X$, c'est le système local $V^\nabla$ des sections horizontales de $\nabla$ (théorème de Cauchy-Kowalewski et lemme de Poincaré holomorphe); (\ref{enum:caslisse}$'$) si de plus $(V,\nabla)$ sous-tend \emph{une variation de $\QQ$-structure de Hodge polarisable}, le théorème de Lefschetz difficile $L^k:H^{\dim X-k}(X,V^\nabla)\isom H^{\dim X+k}(X,V^\nabla)$ a été montré par Deligne (voir \cite[Th\ptbl2.9]{Zucker79}), le théorème de Lefschetz difficile proprement dit, démontré par Hodge (voir \cite{Hodge50}) étant le cas $(V,\nabla)=(\cO_X,d)$.
\item\label{enum:casouvert}
L'extension de ce résultat au cas où $(V,\nabla)$ est un fibré holomorphe à connexion intégrable sur le complémentaire $X^o$ d'une hypersurface $D$ de $X$ et satisfait (\ref{enum:caslisse}$'$) a fait l'objet de nombreux travaux (\cite{Schmid73,Zucker79,C-K82,C-K-S86,C-K-S87, Kashiwara85, K-K87}), aboutissant au théorème de Lefschetz difficile pour la cohomologie d'intersection sur~$X$ du système local $V^\nabla$, lorsque $D=X\moins X^o$ est un diviseur à croisements normaux.
\item\label{enum:cassing}
M\ptbl Saito \cite{MSaito86} a supprimé l'hypothèse~$X$ lisse et~$D$ à croisements normaux en introduisant la catégorie des $\cD$\nobreakdash-modules de Hodge polarisables. Le théorème \ref{th:HLT} s'applique aux complexes $\DR\cM=\IC_Z(\cL)[\dim Z]$ pourvu que $\cL$ sous-tende une variation de $\QQ$-structure de Hodge polarisable (voir aussi \cite{B-B-D81,C-M03,C-M09} pour d'autres approches dans le cas des systèmes locaux d'origine géométrique).
\end{enumeratea}

Dans \eqref{enum:casouvert} et \eqref{enum:cassing}, on travaille avec le prolongement de Deligne de $(V,\nabla)$, qui est singularité régulière à l'infini (\ie le long de~$D$), et ce comportement \og modéré\fg est nécessaire \emph{a priori} pour appliquer les méthodes de la théorie de Hodge, d'après le théorème de régularité de Griffiths-Schmid \cite[Th\ptbl4.13]{Schmid73}). Par contre, la notion de \emph{variation de structure de twisteur polarisée}, introduite par Simpson \cite{Simpson97} autorise des singularités irrégulières à l'infini, et permet d'aborder, par une théorie de Hodge \og sauvage\fg, le théorème \ref{th:HLT} pour les $\cD$-modules holonomes simples à singularités éventuellement irrégulières.

Cette notion de variation de structure de twisteur polarisée intervient déjà en l'absence de singularité (sous la forme d'une \emph{métrique harmonique}, voir le dictionnaire du~\S\ref{sec:dictionnaire}). Pour une variation de structure de Hodge polarisée, c'est la structure qui reste lorsqu'on ne garde que la connexion plate et la métrique hermitienne de polarisation. Sous la seule hypothèse de semi-simplicité de $(V,\nabla)$ ou, de manière équivalente, du faisceau localement constant $V^\nabla$, le cas lisse \eqref{enum:caslisse} du théorème \ref{th:HLT} provient de l'existence, due à Corlette \cite{Corlette88}, d'une métrique dite \emph{harmonique} pour $(V,\nabla)$ (le cas où $(V,\nabla)$ est unitaire ou, plus généralement, une variation de structure de Hodge polarisée, en étant un cas très particulier). On peut en effet développer dans ce cadre la théorie harmonique du laplacien et obtenir les identités de Kähler, qui conduisent au théorème~\ref{th:HLT} (voir \cite{Simpson92}). Il faut noter que l'hypothèse de semi-simplicité de $V^\nabla$ est importante, et il est facile de donner un exemple où \ref{th:HLT} est en défaut sans cette hypothèse: sur une courbe de genre $g\geq2$, toute extension non triviale $\cL$ du faisceau constant $\CC$ par un système local non constant de rang $1$ satisfait à $\dim H^0(X,\cL)\neq\dim H^2(X,\cL)$.

Nonobstant la remarque \ref{rem:H0}, Simpson \cite{Simpson90} a montré, dans le cas où~$X$ est une courbe, l'existence d'une métrique harmonique~$h$ pour $(V,\nabla)$ sur $X^o\subset X$, avec un comportement modéré aux points de~$D$ (voir aussi \cite{Biquard91,J-Z97} en dimension $\geq1$). L'analyse asymptotique qu'il fait de cette métrique au voisinage de~$D$ prolonge à ce cadre celle faite par Schmid \cite{Schmid73} dans le cas des variations de structures de Hodge polarisées, ce qui permet notamment de calculer la cohomologie d'intersection $H^1(X,j_*V^\nabla)$ ($j:X^o\hto X$) comme un espace de cohomologie $L^2$ relativement à~$h$ et à une métrique de type Poincaré sur~$X^o$, et qui prolonge les résultats de Zucker \cite{Zucker79} (\cite{Biquard97}, voir aussi \cite[\S6.2]{Bibi01c}, \cite[\S20.2]{Mochizuki07}, \cite{J-Y-Z07}). Ceci aboutit à un énoncé analogue à la dégénérescence en $E_1$ de la suite spectrale Hodge $\implique$ de~Rham, à savoir le calcul de $\dim H^1(X,j_*V^\nabla)$ en terme de la cohomologie de Dolbeault du fibré de Higgs parabolique associé.

Le cas \og modéré\fg du théorème \ref{th:HLT} est celui où le $\cD_X$-module $\cM$ est à singularités régulières, \ie $\DR\cM=\IC_Z(\cL)[\dim Z]$ avec~$Z$ irréductible et $\cL$ simple sur~$Z^o$. T\ptbl Mochizuki a résolu ce cas dans \cite{Mochizuki07} en étendant les méthodes évoquées ci-dessus. La stratégie de la démonstration, qui vaut aussi pour le cas \og sauvage\fg, c'est-à-dire lorsque $\cM$ est à singularités irrégulières, sera expliquée au~\S\ref{sec:strategie}. Le cas modéré a aussi été résolu par Drinfeld \cite{Drinfeld01} par une méthode de réduction à la caractéristique~$p$ réminiscente de \cite{B-B-D81}. Drinfeld s'appuyait cependant sur une conjecture faite par de~Jong \cite{Jong01}, démontrée depuis \cite{B-K03,Gaitsgory04}.

Récemment, Krämer et Weissauer \cite{K-W11} ont utilisé le cas modéré de \ref{th:HLT} pour montrer un théorème d'annulation, pour tout faisceau pervers $\cF$ sur une variété abélienne complexe~$X$, des espaces $\bH^j(X,\cF\otimes\cL)$ pour tout $j\neq0$ et presque tout système local~$\cL$ de rang $1$. 

La suite du texte insistera donc sur les nouveaux outils introduits  par T\ptbl Mochizuki \cite{Mochizuki08} (après ceux de \cite{Mochizuki07}) pour passer du cas \og modéré\fg au cas \og sauvage\fg.

\subsubsection*{Remerciements}
Ils vont à T\ptbl Mochizuki, ainsi qu'à M\ptbl Saito, Ch\ptbl Schnell, Ch\ptbl Seven\-heck et C\ptbl Simpson pour les multiples suggestions qui m'ont aidé à améliorer la version préliminaire de ce texte.

\section{Dictionnaire}\label{sec:dictionnaire}
Je vais expliciter le dictionnaire fibré plat harmonique/fibré de Higgs harmonique/variation de structure de twisteur polarisée pure de poids~$0$ (voir \cite{Simpson92,Simpson97}), car il est essentiel pour le cas singulier. Les équations \eqref{eq:pseudocourbure} ci-dessous et l'idée de la construction sur l'espace twistoriel d'un fibré avec $\hb$-connexion remontent à Hitchin \cite{Hitchin87}. Drinfeld m'a aussi indiqué les travaux de Zakharov, Mikhailov et Shabat \cite{Z-M78,Z-S79} où on trouve ce type d'équation sous le nom de \og chiral field equations\fg.

Soient $(V,\nabla)$ un fibré holomorphe muni d'une connexion holomorphe intégrable, et $(H,D=\nabla+\ov\partial)$ le fibré $C^\infty$ plat associé. À toute métrique hermitienne~$h$ sur $H$ on associe (\cf\cite{Simpson92}) une unique connexion métrique $\partial_E+\ov\partial_E$ caractérisée par le fait que, si on considère les deux morphismes $\cC^\infty_X$-linéaires $\theta:=\nabla-\partial_E:H\to\cA_X^{1,0}\otimes H$ et $\theta^\dag:=\ov\partial-\ov\partial_E:H\to\cA_X^{0,1}\otimes H$, alors $\theta^\dag$ est le~$h$-adjoint de $\theta$ (si~$D$ est déjà compatible à~$h$ on a $\theta=0$, $\theta^\dag=0$). On dit que $(V,\nabla,h)$ est un \emph{fibré plat harmonique}\footnote{\label{footnote:harmonique}Nous suivons ici la terminologie de \cite{Simpson92}; Mochizuki emploie le terme \og pluri-harmonique\fg pour l'équation \eqref{eq:pseudocourbure}, afin de la distinguer de l'équation \emph{a priori} plus faible sur une variété kählérienne, aussi considérée dans \cite{Corlette88,Simpson92}, à savoir $\Lambda G(\nabla,h)=0$.} si sa \emph{pseudo-courbure} $G(\nabla,h)$ est nulle:
\begin{equation}\label{eq:pseudocourbure}
G(\nabla,h):=-4(\ov\partial_E+\theta)^2=0,\quad\text{\ie\ }\ov\partial_E^2=0,\quad\ov\partial_E(\theta)=0,\quad\theta\wedge\theta=0.
\end{equation}
(Ces trois conditions sont redondantes, les deux dernières impliquant la première; sur une variété kählérienne compacte, on peut même se contenter de la seconde, voir \cite[Rem\ptbl21.33 \& Prop\ptbl21.39]{Mochizuki07}.) Pour un fibré plat harmonique, $E\defin\ker\ov\partial_E:H\to\cA_X^{0,1}\otimes H$ est un fibré holomorphe, et $\theta:E\to\Omega^1_X\otimes E$ est un morphisme holomorphe, qui satisfait à $\theta\wedge\theta=0$. Ainsi, $(E,\theta)$ est un \emph{fibré de Higgs} holomorphe.

Partant maintenant d'un fibré de Higgs holomorphe $(E,\theta)$ (\ie $\theta\wedge\theta=0$) et d'une métrique hermitienne~$h$ sur $E$, on dit que $(E,\theta,h)$ est un \emph{fibré de Higgs harmonique} si, notant $\partial_E+\ov\partial_E$ la connexion de Chern associée à la métrique~$h$ sur le fibré holomorphe~$E$, et $\theta^\dag$ le~$h$-adjoint de $\theta$, alors la connexion $\partial_E+\ov\partial_E+\theta+\theta^\dag$ sur $H:=\cC^\infty_X\otimes_{\cO_X}E$ est \emph{intégrable}.

On a ainsi une correspondance bi-univoque
\[
\text{fibré plat harmonique}\longleftrightarrow \text{fibré de Higgs harmonique.}
\]

Lorsque $(V,\nabla,h)$ ou $(E,\theta,h)$ sont harmoniques, il y a en fait une famille à un paramètre de fibrés holomorphes plats qui dégénère sur le fibré de Higgs associé: pour tout $\hb\in\CC$, on pose $V^\hb=\ker(\ov\partial_E+\hb\theta^\dag:H\to\cA_X^{0,1}\otimes H)$, muni de l'opérateur $\nabla^\hb:=\hb\partial_E+\theta$, qu'on appelle une \emph{$\hb$-connexion}. Si $\hb\neq0$, l'opérateur $\frac1\hb\nabla^\hb$ est une connexion holomorphe intégrable sur $V^\hb$ tandis que, si $\hb=0$, on retrouve le fibré de Higgs $(E,\theta)$.

\begin{exem}[le cas de rang $1$ sur un disque épointé]\label{exem:rankone}
Considérons le cas d'un fibré (trivial) de rang $1$ sur le disque unité épointé $\Delta^*$ de coordonnée $z$, muni d'une métrique hermitienne~$h$. Nous noterons $\gU$ l'ensemble des classes d'équivalence de couples $\gu=(a,\alpha)\in\RR\times\CC$ modulo $\ZZ\times\{0\}$. Alors:
\par\smallskip
\emph{L'ensemble des classes d'isomorphisme de fibrés de Higgs (ou de fibrés plats) harmoniques de rang $1$ sur $\Delta^*$ est en correspondance bijective avec l'ensemble des couples $(\psi,\gu)$, avec $\psi\in\cO(\Delta^*)$ sans terme constant et $(\gu\bmod\ZZ)\in\gU$.}

\begin{proof}
Nous la ferons dans le cas Higgs, le cas plat étant similaire. Soit $(E,\theta,h)$ un fibré de Higgs harmonique de rang $1$ sur $\Delta^*$. Nous allons lui associer un unique couple $(\psi,\gu\bmod\ZZ)$. Soit $\epsilon$ une base holomorphe de $E$. On a
\[
\theta \epsilon=\varphi(z)\,\epsilon\,\rd z,\quad \text{$\varphi(z)$ holomorphe sur $\Delta^*$.}
\]
Posons $\varphi(z)=\partial_z\psi(z)+\alpha/z$ avec $\psi\in\cO(\Delta^*)$ sans terme constant et $\alpha\in\CC$. Posons aussi $\norme{\epsilon}_h=\exp(\eta(z))$ où $\eta$ est réelle et $C^\infty$ sur $\Delta^*$. On peut vérifier que la condition d'harmonicité de $(E,\theta,h)$ équivaut au fait que la fonction $\eta$ est \emph{harmonique} sur $\Delta^*$. Elle s'écrit donc $\reel\gamma(z)-a\log|z|$ avec $\gamma$ holomorphe sur $\Delta^*$ et $a\in\RR$. Remplaçant~$\epsilon$ par $e=\exp(-\gamma(z))\cdot \epsilon$, on peut supposer que $\eta(z)=-a\log|z|$ avec $a\in\RR$, et on a alors $\norme{e}_h=|z|^{-a}$. On a ainsi obtenu un couple $(\psi,\gu)$.

On calcule ensuite que $v:=|z|^{-2\ov\alpha}\exp(\psi-\ov\psi)\cdot e$ est une base holomorphe du fibré plat $V$ associé, de norme $\norme{v}_h=|z|^{-a-2\reel\alpha}$. De plus,
\[
\theta e=(z\partial_z\psi+\alpha)\frac{\rd z}{z}\otimes e,\quad \nabla v=(2z\partial_z\psi+\ge(1,\gu))\frac{\rd z}{z}\otimes v,\quad\ge(1,\gu):=-a+2i\im\alpha.
\]
Soient $\epsilon'$ une autre base holomorphe de $E$, et $e'$ construite comme plus haut avec~$\Vert e'\Vert_h=|z|^{-a'}$ pour un certain \hbox{$a'\in\RR$}. D'où un couple $(\psi',\gu')$. Alors $e'=\nu(z)e$ avec $\nu(z)$ holomorphe et à croissance modérée, donc méromorphe, d'où $a'-a\in\ZZ$. Le champ de Higgs a la même expression dans les bases~$e$ et $e'$, ce qui implique $\psi=\psi'$ et $\alpha=\alpha'$.
\end{proof}

Posons maintenant, pour tout $\hb\in\CC$ fixé,
\bgroup\numstareq
\begin{equation}\label{eq:pe}
\gp(\hb,\gu)=a+2\reel(\ov\alpha\hb),\quad\ge(\hb,\gu)=\alpha-a\hb-\ov\alpha\hb^2,\quad\text{et}\quad v^\hb= e^{\ov\hb\psi-\hb\ov\psi}|z|^{-2\ov\alpha\hb}\cdot e.
\end{equation}
\egroup
Alors $v^\hb$ est une base holomorphe de $V^\hb$, de norme $\Vert v^\hb\Vert_h=|z|^{-\gp(\hb,\gu)}$ et
\bgroup\numstarstareq
\begin{equation}\label{eq:vhb}
\nabla^\hb v^\hb=((1+|\hb|)^2z\partial_z\psi+\ge(\hb,\gu))\frac{\rd z}{z}\otimes v^\hb.
\end{equation}
\egroup
\end{exem}

Revenons à la situation générale. Ces deux notions équivalentes (fibré plat harmonique et fibré de Higgs harmonique) sont aussi équivalentes à la notion de \emph{variation de structure de twisteur polarisée pure de poids~$0$}. Pour la définir, introduisons la droite projective $\PP^1$ munie de deux cartes affines $\CC_\hb,\CC_\mu$ de coordonnées respectives $\hb$ et $\mu$, avec $\mu=1/\hb$ sur l'intersection des deux cartes. La présentation qui suit n'est pas exactement celle donnée par Simpson \cite{Simpson97}, mais lui est équivalente et sera plus commode dans les situations singulières. Il s'agit de décrire des fibrés sur $X\times\PP^1$ qui sont holomorphes par rapport à $\PP^1$ et $C^\infty$ par rapport à $X$. Cette présentation permet de ne travailler qu'avec des fibrés holomorphes sur $X\times\CC_\hb$ et, plus loin (\S\ref{sec:Dmodtw}), avec des $\cR_{X\times\CC_\hb}$-modules holonomes.

Soit $\sigma:\PP^1\to\ov\PP^1$ l'involution anti-holomorphe qui prolonge continûment l'application $\CC_\mu\to\CC_\hb$, $\mu\mto\hb=-1/\ov\mu$. Si $f(x,\hb)$ est holomorphe en~$x$ et $\hb$, alors la fonction $(\sigma^*\ov f)(x,\mu)=\ov{f(x,-1/\ov\mu)}$ est anti-holomorphe en~$x$ et holomorphe en~$\mu$. Si $\cH$ est un fibré holomorphe sur $X\times\CC_\hb$, alors $\sigma^*\ov\cH$ est un fibré holomorphe sur $\ov X\times\CC_\mu$, où~$\ov X$ est la variété complexe conjuguée de~$X$. Si $\cH',\cH''$ sont deux fibrés holomorphes sur $X\times\CC_\hb$, un \emph{pré-recollement} (entre le dual~$\cH^{\prime\vee}$ et $\sigma^*\ov\cH{}''$) est un accouplement $\cO_{X\times\bS}\otimes_{\cO_{\bS}}\cO_{\ov X\times\bS}$-linéaire
\[
C:\cH'_{|X\times\bS}\otimes_{\cO_{\bS}}\sigma^*\ov\cH{}''_{|X\times\bS}\to\cC^{\infty,\an}_{X\times\bS},
\]
où $\bS=\{|\hb|=1\}$, $\cO_{\bS}$ est le germe le long de $\bS$ de $\cO_{\CC_\hb}$ et $\cC^{\infty,\an}_{X\times\bS}$ le faisceau des germes le long de $X\times\bS$ de fonctions $C^\infty$ holomorphes par rapport à $\hb$. Si $\cH'$ et~$\cH''$ sont munis de $\hb$\nobreakdash-connexions, on demande que $C$ soit compatible en un sens naturel. Ces triplets forment naturellement une catégorie. L'\emph{adjoint} $(\cH',\cH'',C)^*$ est le triplet $(\cH'',\cH',C^*)$, avec $C^*(u'',\sigma^*\ov u{}')\defin \ov{C(u',\sigma^*\ov u{}'')}$, et les $\hb$-connexions restent compatibles à $C^*$.

Une \emph{variation de structure de twisteur} est une telle donnée $(\cH',\cH'',C)$ avec $\hb$\nobreakdash-connexions intégrables telle que, pour tout $x\in X$, l'accouplement induit par $C$ soit non dégénéré et définisse donc un fibré holomorphe sur $\PP^1$ par recollement de $\cH^{\prime\vee}_{|\{x\}\times\CC_\hb}$ et $\sigma^*\ov\cH{}''_{|\{x\}\times\CC_\mu}$. C'est une variation de structure de twisteur \emph{pure de poids $w\in\ZZ$} si, pour tout~$x$, le fibré obtenu est isomorphe à une puissance de $\cO_{\PP^1}(w)$. Si le poids $w$ est nul, les sections globales de ce fibré à~$x$ fixé forment un espace vectoriel~$H_x$ de dimension égale au rang de $\cH'$ et $\cH''$, et l'adjoint a pour sections globales l'espace adjoint~$\ov H{}_x^\vee$. On définit alors une \emph{polarisation} comme un isomorphisme~$\cS$ de $(\cH',\cH'',C)$ sur son adjoint $(\cH',\cH'',C)^*$, compatible aux $\hb$-connexions, tel que, pour tout~$x$, l'isomorphisme induit $H_x\isom \ov H{}_x^\vee$, vu comme un accouplement sesquilinéaire sur $H_x$, soit une forme hermitienne \emph{définie positive}. Le fibré $H$ sur~$X$ dont les fibres sont les $H_x$ est alors un fibré $C^\infty$ muni d'une métrique hermitienne~$h$.

\begin{lemm}[\cite{Simpson97}]
Soit $(\cH',\cH'',C,\cS)$ une variation de structure de twisteur polarisée de poids~$0$. La restriction $\cH''$ à $\hb=1$ (\resp $\hb=0$) munie de la connexion (\resp le champ de Higgs) induite par la $\hb$-connexion est un fibré holomorphe plat (\resp de Higgs) de fibré $C^\infty$ sous-jacent isomorphe à $H$, et la métrique~$h$ en fait un fibré plat (\resp de Higgs) harmonique.

Réciproquement, la construction $V^\lambda$ à partir d'un fibré plat (\resp de Higgs) harmonique permet de définir une variation de structure de twisteur polarisée pure de poids~$0$ en posant $\cH'=\cH''=\ker(\ov\partial_\hb+\ov\partial_E+\hb\theta^\dag:\cC^\infty_{X\times\CC_\hb}\otimes_{\cC^\infty_X}H\to\cA_{X\times\CC_\hb}^{0,1}\otimes_{\cC^\infty_X} H)$, $\cS=\id$, $C$ est induit naturellement par~$h$ et la $\hb$-connexion par $\hb\partial_E+\theta$.
\end{lemm}

\begin{enonce*}[remark]{Exemple \ref{exem:rankone}, suite}
Considérons maintenant $v^\hb$ comme dépendant de $\hb$. Alors la formule \eqref{eq:vhb} montre que $\nabla^\hb$ ne s'exprime de manière holomorphe en $\hb$ que si $\psi=0$. Si $\psi\neq0$, on peut considérer la base holomorphe $\wt v{}^\hb=e^{-|\hb|^2\psi}\cdot v^\hb$ pour corriger le problème. On constate alors que, d'une part, il n'y a pas unicité de choix (on pourrait tout aussi bien prendre $e^{c-|\hb|^2\psi}\cdot v^\hb$, avec $c\in\CC$), et d'autre part tous les choix conduisent à une base dont la norme n'est plus à croissance modérée à l'origine, si $\hb\neq0$ et $\psi$ n'est pas holomorphe à l'origine. Mais, dans ce cas, en faisant dépendre la métrique de $\hb$ de manière convenable, on retrouve la propriété de croissance modérée par rapport à cette métrique modifiée (on retrouvera ceci au point~\eqref{enum:hbvariablee} du \S\ref{subsec:hbvariable}).
\end{enonce*}

\section{Stratégie de la démonstration}\label{sec:strategie}
La stratégie utilisée pour démontrer le théorème \ref{th:HLT} suit celle de M\ptbl Saito \cite{MSaito86}.

\begin{enumerate}
\item\label{enum:programme1}
La première étape consiste à enrichir la structure de $\cD_X$-module holonome (où~$X$ est une variété analytique complexe quelconque) afin de pouvoir disposer d'une notion de poids.
\begin{itemize}
\item
La catégorie des $\cD_X$-modules holonomes munis d'une bonne filtration et d'une structure rationnelle (\ie un isomorphisme du complexe de de~Rham et du complexifié d'un $\QQ$-faisceau pervers) a été considérée par Saito \cite{MSaito86}. L'oubli de la filtration et de la structure rationnelle définit un foncteur d'oubli vers la catégorie des $\cD$-modules holonomes.
\item
Dans le cas présent, on généralise les objets $(\cH',\cH'',C)$ du \S\ref{sec:dictionnaire}: notant $\cR_{X\times\CC}$ le faisceau des opérateurs $\hb$-différentiels engendré par les fonctions $\cO_{X\times\CC}$ et les $\hb$-champs de vecteurs $\hb\partial_{z_i}$, on considère les triplets $(\cM',\cM'',C)$, où $\cM',\cM''$ sont des $\cR_{X\times\bS}$-modules holonomes (en un sens naturel) et $C$ est un accouplement entre $\cM'_{|X\times\bS}$ et $\sigma^*\ov\cM{}''_{|X\times\bS}$ à valeurs dans le faisceau des distributions sur $X\times\nobreak\bS$ qui dépendent continûment de $\bS$. La restriction $\cM''/(\hb-1)\cM''$ définit un foncteur d'oubli à valeurs dans la catégorie des $\cD$-modules holonomes.
\end{itemize}
\item
Sans plus de contraintes, les catégories ci-dessus ne sont pas abéliennes.
\begin{itemize}
\item\label{enum:programme2}
L'idée de Saito consiste essentiellement à imposer des conditions supplémentaires locales: pour tout germe de fonction holomorphe $f$ sur~$X$ le foncteur $\psi_f^\rmod$ des \emph{cycles proches modérés} le long de $f=0$, défini \emph{a priori} sur la catégorie des $\cD$\nobreakdash-modules holonomes à l'aide de la $V$-filtration de Kashiwara-Malgrange, doit exister pour les objets filtrés considérés. Il impose donc d'une part l'existence d'un tel foncteur et d'autre part que le résultat donne un objet du même type (à graduation près par la filtration dite \og monodromique\fg) avec une dimension de support strictement plus petite. Lorsque la dimension du support est nulle, on impose d'obtenir une structure de Hodge polarisée. Le cas le plus simple de ce procédé est la restriction à un point d'une variation de structure de Hodge. On obtient ainsi la catégorie des $\cD$\nobreakdash-mo\-dules de Hodge polarisables \cite{MSaito86}.
\item
Cette idée se transporte assez directement aux triplets $(\cM',\cM'',C)$, la définition de $\psi_f^\rmod$ sur l'accouplement $C$ s'obtenant en prenant le résidu en différentes valeurs de $s$ du transformé de Mellin de la distribution $|f|^{2s}C$. On obtient ainsi la catégorie des $\cD$-modules holonomes avec structure de twisteur polarisable modérée \cite{Bibi01c,Mochizuki07}.
\item
Pour le cas sauvage (singularités irrégulières), l'utilisation des cycles proches modérés est insuffisante. Les cycles proches irréguliers, tels que définis par Deligne \cite{Deligne83} sont par contre suffisants (\cite{Bibi06b}, \cite{Mochizuki08}). On obtient la catégorie des $\cD$-modules holonomes avec structure de twisteur polarisable sauvage \cite{Mochizuki08} (voir le \S\ref{sec:Dmodtw}).
\end{itemize}
\item\label{enum:programme3}
On montre alors le théorème \ref{th:HLT} pour les $\cD$-modules de Hodge polarisables ou avec structure de twisteur polarisable. Les espaces d'hypercohomologie sont soit des espaces complexes filtrés avec une structure rationnelle et une polarisation, soit des triplets $(\cH',\cH'',C)$ formés de $\cO_{\CC_\hb}$-modules et d'un accouplement. On montre de plus que les morphismes de Lefschetz sont strictement compatibles aux filtrations ou restent des isomorphismes par restriction à $\hb=1$. Ceci est obtenu en montrant que les espaces d'hypercohomologie sont des structures de Hodge pures ou des structures de twisteur pures. On en conclut que tous les $\cD$-modules obtenus par oubli de la filtration à partir des $\cD$-modules de Hodge polarisables ou par restriction à $\hb=1$ de $\cD$-modules avec structure de twisteur polarisable satisfont au théorème de Lefschetz difficile.

\item\label{enum:programme4}
Les résultats du \eqref{enum:programme3} sont obtenus par la méthode de pinceaux de Lefschetz. On se ramène ainsi à les montrer dans le cas des courbes.
\begin{itemize}
\item
Pour les $\cD$-modules de Hodge polarisables \cite{MSaito86}, on s'appuie sur les théorèmes calculant la cohomologie $L^2$ obtenus par Zucker \cite{Zucker79}, car la restriction d'un tel $\cD$-module filtré à un ouvert de Zariski dense de la courbe n'est autre qu'une variation de $\QQ$-structure de Hodge polarisable.
\item
Pour les $\cD$-modules avec structure de twisteur polarisable, on commence par remarquer que la restriction à un ouvert de Zariski dense fournit, par le dictionnaire du \S\ref{sec:dictionnaire}, un fibré holomorphe plat avec métrique harmonique. De plus, cette métrique est modérée, au sens de \cite{Simpson90}, ou sauvage, au sens expliqué au \S\ref{sec:sauvage}. Les résultats de \cite{Simpson90} puis de \cite{Mochizuki08} dans le cas des courbes pour de telles métriques permettent d'adapter la méthode de Zucker, comme indiqué dans l'introduction (une approche un peu différente dans le cas modéré est utilisée dans \cite{Biquard97} et \cite{Bibi01c}; voir aussi \cite{Bibi98} pour un lemme de Poincaré $L^2$ dans le cas sauvage).
\end{itemize}

\item\label{enum:programme5}
On arrive maintenant à la deuxième partie du programme, qui est plus analytique. Il s'agit d'identifier exactement les $\cD$-modules holonomes produits au point~\eqref{enum:programme3}. Dans la mesure où on a mis en place, durant la preuve du théorème \ref{th:HLT} pour les catégories de $\cD$-modules de Hodge polarisables ou avec structure de twisteur polarisable, les outils tels que l'image directe par un morphisme projectif et le théorème de décomposition analogue à celui de \cite{B-B-D81}, on peut ramener cette identification au cas où le support~$Z$ est lisse et où l'ouvert de Zariski~$Z^o$ de lissité du $\cD$-module holonome a pour complémentaire un diviseur à croisements normaux.
\begin{itemize}
\item
Pour les $\cD$-modules de Hodge polarisables, M\ptbl Saito \cite[Th\ptbl3.21]{MSaito87} les identifie, via la correspondance de Riemann-Hilbert, à ceux qui correspondent aux complexes d'intersection de variations de structure de Hodge polarisable. Un point délicat est la reconstruction, à partir d'une telle variation, d'un $\cD$-module de Hodge polarisable, et l'ingrédient essentiel est le théorème d'existence d'une structure de Hodge mixte limite, dû à Cattani, Kaplan et Schmid d'une part, Kashiwara et Kawai de l'autre, ainsi que la description de celle-ci, généralisant en toute dimension le théorème de Schmid en dimension~$1$ (voir l'introduction).
\item
Pour les $\cD$-modules avec structure de twisteur polarisable \emph{modérée}, Mochizuki identifie les $\cD$-modules holonomes produits au point \eqref{enum:programme3} à ceux correspondant, par Riemann-Hilbert, aux complexes d'intersection de fibrés plats munis d'une métrique harmonique modérée (voir le \S\ref{sec:sauvage} ci-dessous). La reconstruction (voir le~\S\ref{sec:prolong}) s'appuie aussi, \emph{in fine}, sur la théorie asymptotique des variations de structure de Hodge polarisable.
\item
Dans le cas sauvage, on se ramène d'abord à considérer une situation mieux contrôlée (cas sauvage et bon, voir \S\ref{sec:sauvage}). Mochizuki a montré, comme conséquence de sa démonstration de la correspondance de Hitchin-Kobayashi sauvage (\S\ref{sec:HK}), l'existence d'une compactification convenable $(Z,D)$ de~$Z^o$ pour laquelle cette propriété est satisfaite. On peut désormais déconnecter cet argument de l'ensemble de la démonstration, et utiliser ici les résultats de Kedlaya \cite{Kedlaya09,Kedlaya10} (voir le \S\ref{subsec:singirreg}), ce que je ferai pour simplifier l'exposé.
\end{itemize}

\item\label{enum:programme6}
La correspondance de Hitchin-Kobayashi, expliquée dans ce cadre à la section~\ref{sec:HK}, permet enfin à Mochizuki de mettre en correspondance bijective les fibrés plats harmoniques modérés ou sauvages et bons sur $(Z,D)$ avec les fibrés \emph{méromorphes plats semi-simples} lisses sur~$Z^o$, généralisant ainsi le théorème de Corlette \cite{Corlette88} du cas projectif (voir aussi \cite{J-Z97} pour le cas quasi-projectif). Ceux-ci sont eux-mêmes en correspondance bijective avec les $\cD_Z$-modules holonomes semi-simples lisses sur~$Z^o$. La démonstration aux \oldS\S\ref{sec:HK} et~\ref{sec:prolong} insistera sur l'aspect \og reconstruction\fg (essentielle surjectivité dans le théorème \ref{th:equiv} ci-dessous). Il faut noter cependant que, dans le cas sauvage, l'aspect direct (les $\cD$-modules obtenus au point~\eqref{enum:programme3} sont semi-simples) et la pleine fidélité sont des points non triviaux car, notamment, on ne dispose pas du prolongement méromorphe canonique de Deligne (\cf \S\ref{subsec:singirreg}). Cette question est traitée au \S19.3 de \cite{Mochizuki08}.

En conclusion, les $\cD$-modules holonomes produits au point \eqref{enum:programme3} sont exactement les $\cD$-modules holonomes semi-simples sur~$X$, ce qui termine la démonstration du théorème~\ref{th:HLT}.
\end{enumerate}

\begin{rema}
On a ici perdu un peu de symétrie entre fibrés de Higgs et fibrés plats. En effet, les fibrés de Higgs harmoniques sauvages et bons sur $(Z,D)$ sont en correspondance bijective avec les fibrés de Higgs poly-stables avec nombres caractéristiques paraboliques nuls (voir le \S\ref{subsec:locab}). Mais on ne sait pas identifier les objets de Higgs obtenus au point~\eqref{enum:programme3} par restriction à $\hb=0$, en tant que $\cO_{T^*X}$-modules cohérents.
\end{rema}

\section{Fibrés de Higgs sauvages et fibrés méromorphes plats à singularités irrégulières}\label{sec:sauvage}

\setcounter{subsection}{-1}
\subsection{Convention}\label{subsec:convention}
Dans tout ce texte, nous appellerons \og situation globale\fg la donnée d'une variété \hbox{projective} lisse~$X$ et d'un diviseur~$D$ à croisements normaux dans~$X$. On note \hbox{$j:X^o:=X\moins D\hto X$} l'inclusion. Les composantes du diviseur, supposées lisses, sont indexées par un ensemble fini $\ccI$. On se donne aussi un fibré en droites ample $L$ sur~$X$.

Par \og situation locale\fg, nous entendrons plutôt que~$X$ est un produit $\Delta^n$ de disques muni de coordonnées $(z_1,\dots,z_n)$, et~$D$ a pour équation $z_1\cdots z_\ell=\nobreak0$, de sorte que l'ensemble fini $\ccI$ est ici égal à $\{1,\dots,\ell\}$.

\subsection{Fibrés de Higgs sauvages}\label{subsec:Higgssauvage}
Considérons la situation locale. Soit $(E,\theta)$ un fibré de Higgs holomorphe sur $X^o$. On écrit le champ de Higgs sous la forme $\theta=\sum_{i=1}^\ell F_i\rd z_i/z_i+\sum_{j=\ell+1}^nG_j\rd z_j$, où les $F_i,G_j$ sont des endomorphismes holomorphes de $E$. Les coefficients $f_{i,k},g_{j,k}$ des polynômes caractéristiques des $F_i$ et des $G_j$ sont des fonctions holomorphes sur $X^o$.

Le fibré de Higgs $(E,\theta)$ est \emph{modéré} (dans la carte considérée) si $f_{i,k},g_{j,k}$ se prolongent en des fonctions \emph{holomorphes} sur~$X$ (on notera de la même manière ces prolongements) et $f_{i,k|D_i}$ est \emph{constante} pour tous $i,k$. En particulier, les valeurs propres (et leur multiplicité) de $F_{i|D_i}$ sont constantes et égales à celles de $F_i(0)$. Du fait de la condition de Higgs, les endomorphismes $F_i,G_j$ commutent. On en déduit que $(E,\theta)$ se décompose localement suivant l'ensemble $\Sp(\theta)$ des valeurs propres $(\alpha_1,\dots,\alpha_\ell)$ de $(F_1(0),\dots,F_\ell(0))$ en somme directe de sous-fibrés de Higgs modérés (voir \cite[\S8.2.1]{Mochizuki07}):
\refstepcounter{defi}
\bgroup\numstareq
\begin{equation}
\label{eq:decompvp}
(E,\theta)\simeq\bigoplus_{\alphag\in\Sp(\theta)}(E_\alphag,\theta_\alphag).
\end{equation}
\egroup

Dans le cas \og sauvage\fg, les fonctions $f_{i,k},g_{j,k}$ sont autorisées à avoir des pôles le long de~$D$, mais de manière contrôlée. Travaillant toujours dans des coordonnées locales, on notera $\cO(*D)$ l'espace des fonctions méromorphes sur~$X$ à pôles d'ordre arbitraire le long de~$D$, et on s'intéressera aux parties polaires $\cO(*D)/\cO$.

Le fibré de Higgs $(E,\theta)$ a une \emph{décomposition sauvage sans ramification} s'il existe une famille finie $\Irr(\theta)\subset\cO(*D)/\cO$ de parties polaires et une décomposition
\bgroup\numstarstareq
\begin{equation}
\label{eq:decomp}
(E,\theta)=\bigoplus_{\ga\in\Irr(\theta)}(E_\ga,\theta_\ga)
\end{equation}
\egroup
telle que, pour chaque $\ga\in\Irr(\theta)$, le fibré de Higgs $(E_\ga,\theta_\ga-\rd\ga\otimes\id)$ soit \emph{modéré} (cette condition ne dépend que de la partie polaire $\ga$, et pas d'un relèvement à $\cO(*D)$, d'où le raccourci de notation). Le fibré de Higgs $(E,\theta)$ a une \emph{décomposition sauvage avec ramification} s'il admet une décomposition \eqref{eq:decomp} après image inverse par un morphisme fini ramifié autour de~$D$, décrit en coordonnées locales convenables par $\rho:(x_1,\dots,x_n)\mto(z_1,\dots,z_n)=(x_1^{\nu_1},\dots,x_\ell^{\nu_\ell},x_{\ell+1},\dots,x_n)$. Enfin, $(E,\theta)$ est dit \emph{sauvage} si, pour tout point de $D$, il existe une modification projective $\pi$ d'un voisinage de ce point telle que $\pi^{-1}(D)$ soit encore à croisements normaux et $\pi^*(E,\theta)$ admette une décomposition sauvage avec ramification au voisinage de tout point de $\pi^{-1}(D)$.

\begin{enonce*}[remark]{Exemple \ref{exem:rankone}, suite}
Puisque $\psi$ est sans terme constant, on peut écrire $\psi=\rd(\ga+\eta)$ avec $\eta$ holomorphe et $\ga$ holomorphe en $z^{-1}$ sans terme constant. Le fibré de Higgs $(E,\theta)$ est modéré si et seulement si $\ga=0$. Il est sauvage si et seulement si $\psi$ est méromorphe en $z=0$ (\ie $\ga\in z^{-1}\CC[z]$). 
\end{enonce*}

\skpt
\begin{remas}\label{rem:Higgsmodsauv}
\begin{enumeratei}
\item\label{rem:Higgsmodsauv1}
Ce sera une constante du traitement \og modéré/sauvage\fg que de ramener le cas sauvage au cas modéré par adjonction de tels $\rd\ga$ et sommes directes, éventuellement avec ramification. Le cas des fibrés plats se compliquera par l'introduction de structures de Stokes, qui n'apparaissent pas pour les fibrés de Higgs.
\item\label{rem:Higgsmodsauv2}
Les conditions ci-dessus ne dépendent pas du choix des coordonnées adaptées à~$D$, et si $\dim X=1$ elles se réduisent à l'holomorphie (\resp la méromorphie) sur~$X$ des coefficients du polynôme caractéristique de $F_1$. En dimension $\geq2$, la propriété de \og décomposition sauvage avec ramification\fg implique qu'après ramification, le polynôme caractéristique $\chi_{F_i}(T)$ de~$F_i$ se décompose en $\prod_\ga P_{i,\ga}(T-z_i\partial_{z_i}\ga)$, où les $P_{i,\ga}(T)$ sont à coefficients holomorphes pour tous~$i,\ga$, et $\chi_{G_j}(T)=\prod_\ga P_{j,\ga}(T-\partial_{z_i}\ga)$, où les $P_{j,\ga}$ sont à coefficients holomorphes. Cette dernière propriété est \emph{a~priori} plus forte que la méromorphie des coefficients des $\chi_{F_i},\chi_{G_j}$.
\item\label{rem:Higgsmodsauv3}
Les conditions \og modéré\fg ou \og sauvage\fg ne disent rien sur l'éventuel prolongement du fibré $E$ en un fibré sur~$X$ ni, le cas échéant, du prolongement des $F_i,G_j$ comme endomorphismes méromorphes ou holomorphes de ce fibré prolongé.
\item\label{rem:Higgsmodsauv4}
Ces conditions sont préservées par image inverse par un morphisme $f:(X',D')\to(X,D)$ de variétés munies d'un diviseur à croisements normaux  avec $D'=f^{-1}(D)$. Inversement, étant donné $(E,\theta)$ sur $X^o$, \emph{existe-t-il une modification propre $\pi:X'\to X$ qui est un isomorphisme au-dessus de $X^o$, telle que $X'\moins X^o$ soit un diviseur à croisements normaux, et que $(E,\theta)$ admette une décomposition sauvage avec ramification sur $X'$?} Puisqu'une telle modification est un isomorphisme hors d'un ensemble de codimension $1$ dans $D$, il est nécessaire \emph{a priori} que $\theta$ soit sauvage au voisinage de tout point d'un ouvert de Zariski dense de $D$. Réciproquement, Mochizuki montre \cite[Chap\ptbl15]{Mochizuki08} que c'est bien le cas si les propriétés suivantes sont satisfaites:
\begin{itemize}
\item
$D$ est à croisements normaux dans~$X$,
\item
$\theta$ admet une décomposition sauvage avec ramification génériquement le long de~$D$, 
\item
au voisinage de tout point de~$D$, et dans des coordonnées locales adaptées, on~a, après ramification finie éventuelle, une décomposition $\chi_{F_i}(T)=\prod_{\ga\in A_i} P_{i,\ga}(T-\nobreak z_i\partial_{z_i}\ga)$ comme plus haut.
\end{itemize}
\end{enumeratei}
\end{remas}

\subsection{Fibrés plats à singularités irrégulières}\label{subsec:singirreg}
Arrêtons-nous un instant pour imaginer l'analogue des deux propriétés \og modéré/\allowbreak sauvage\fg pour un fibré holomorphe plat $(V,\nabla)$ sur $X^o$ lorsque~$D$ est à croisements normaux. Considérer comme ci-dessus les polynômes caractéristiques des coefficients de la matrice de connexion n'a plus de sens, mais il existe un unique prolongement méromorphe\footnote{\ie un $\cO_X(*D)$-module localement libre de rang fini $\cV$ muni d'une connexion intégrable $\nabla$.} de $(V,\nabla)$ dans une base duquel les sections horizontales de la connexion ont des coefficients à croissance modérée (prolongement de Deligne). De plus, il existe un prolongement holomorphe canonique sur lequel la connexion est à pôles logarithmiques, et les résidus le long des composantes de~$D$ sont constants. C'est la situation \og modérée\fg (à singularités régulières), qui se comporte donc mieux que l'analogue pour les fibrés de Higgs. L'analogue de la propriété \og sauvage\fg nécessite par contre un prolongement méromorphe de référence du fibré $V$ pour être définie, contrairement au cas Higgs. Tout prolongement méromorphe $(\cV,\nabla)$ distinct du prolongement de Deligne sera dit à \emph{singularités irrégulières}. Dans la situation de l'introduction, si $Z$ est lisse, la condition d'algébricité de $(V,\nabla)$ fournit un prolongement méromorphe bien déterminé du fibré holomorphe associé  $(V,\nabla)^\an$.

La propriété de décomposition (après ramification) paramétrée par des parties polaires $\ga$ analogue à \eqref{eq:decomp} n'est pas satisfaite en général, même si $\dim X=1$. Dans ce cas, elle n'est satisfaite que si on autorise des changements de jauge \emph{formels} pour la connexion.

Pour $\dim X=2$, il était conjecturé (et démontré dans quelques cas particuliers) dans \cite[Conj\ptbl2.5.1]{Bibi97}, qu'une telle propriété est satisfaite après éclatements de~$X$. Cette propriété a été démontrée par T\ptbl Mochizuki \cite{Mochizuki07b} lorsque la connexion est définie algébriquement, par une méthode de réduction à la caractéristique $p$ et une utilisation de la $p$-courbure comme ersatz d'un champ de Higgs, pour laquelle on peut appliquer un énoncé du type de celui de la remarque \ref{rem:Higgsmodsauv}\eqref{rem:Higgsmodsauv4}. Quelque temps après, Kedlaya \cite{Kedlaya09} a proposé une démonstration complètement différente du même énoncé, sans condition d'algébricité de la connexion. Elle repose sur des techniques inspirées des équations différentielles $p$-adiques.

Dans \cite[Th\ptbl16.2.1]{Mochizuki08}, T\ptbl Mochizuki a étendu ce résultat en toute dimension (toujours pour une connexion définie algébriquement), en s'appuyant sur le résultat en dimension $2$ et en utilisant une métrique harmonique pour faire l'aller-retour entre fibré plat et fibré de Higgs par la correspondance de Hitchin-Kobayashi sauvage expliquée plus bas. De ce fait, la démonstration contient une grosse partie d'analyse. Par ailleurs, Kedlaya \cite{Kedlaya10} a aussi pu étendre ses propres méthodes à la dimension quelconque, sans hypothèse d'algébricité.

\begin{theo}[\cite{Mochizuki07b,Mochizuki08}, \cite{Kedlaya09,Kedlaya10}]\label{th:LT}
Soit $X'$ une variété algébri\-que lisse (\resp un germe de variété analytique complexe) et soit $(\cV,\nabla)$ un fibré méromorphe sur $X'$ à connexion intégrable, holomorphe sur un ouvert de Zariski $X^{\prime o}$ de~$X'$. Il existe alors une modification projective $\pi:X\to X'$ avec~$X$ lisse, qui est un isomorphisme au-dessus de $X^{\prime o}$, telle que $X\moins X^{\prime o}$ soit un diviseur à croisements normaux~$D$, et qu'en tout point $x\in D$, le formalisé $(\wh\cO_{X,x}\otimes_{\cO_x}\cV,\wh\nabla)$ se décompose, après ramification éventuelle autour des composantes locales de~$D$, sous la forme
\begin{equation}\tag{$\ref{th:LT}\,*$}\label{eq:LT}
(\wh\cO_{X,x}\otimes_{\cO_x}\cV,\wh\nabla)=\bigoplus_{\ga\in\Irr_x(\nabla)}(\wh\cV_\ga,\wh\nabla_\ga),
\end{equation}
où $\wh\nabla{}_\ga^\reg:=\wh\nabla_\ga-\rd\ga\otimes\id_{\wh\cV_\ga}$ est à singularités régulières (\cite{Deligne70}).
\end{theo}

Ce théorème, tel que démontré par ces auteurs dans sa version plus précise avec la propriété (Bon) ci-dessous, était le chaînon manquant pour analyser les singularités irrégulières de systèmes holonomes d'équations aux dérivées partielles.

\begin{enonce*}[remark]{Exemple \ref{exem:rankone}, suite}
Pour $\hb\neq0$, la dichotomie singularité régulière/irrégulière pour $(V^\hb,\frac1\hb\nabla^\hb)$ est la même que la dichotomie modéré/sauvage du cas Higgs.
\end{enonce*}

\subsection{La condition \og sauvage et bon\fg}\label{subsec:goodwild}

La propriété de décomposition \eqref{eq:decomp} (après ramification locale autour des composantes du diviseur~$D$) pour un fibré de Higgs  ou, en prenant des coefficients formels, pour un fibré méromorphe à connexion plate, est encore insuffisante lorsque $n=\dim X\geq2$ pour l'analyse des propriétés asymptotiques d'une métrique harmonique ou du phénomène de Stokes. Par exemple, dans le cas de deux variables $x_1,x_2$, on cherche à éviter l'existence de sections horizontales de la connexion qui ont un comportement en $\exp(x_1/x_2)$ à cause de la forme \og indéterminée\fg de la limite de $x_1/x_2$ quand $x_1,x_2\to0$. Par contre, on accepte $\exp(1/x_2)$ ou $\exp(1/x_1x_2)$.

Dans des coordonnées locales adaptées au diviseur comme au \S\ref{subsec:Higgssauvage}, on associe à toute partie polaire $\ga\in\cO(*D)/\cO$, écrite sous la forme $\sum_{\bmm\in\ZZ^\ell\times\NN^{n-\ell}}\ga_\bmm z^{\bmm}$, le polyèdre de $\RR^n$ enveloppe convexe des octants $\RR_+^n$ (pour négliger $\cO$) et $\bmm+\RR_+^n$ pour lesquels $\ga_\bmm\neq0$.

Une famille finie $S$ de parties polaires, telle que $\Irr(\theta)$, est dite \emph{bonne} si la propriété suivante est satisfaite:
\enum{.05}{.92}{(Bon)}{les polyèdres de Newton des parties polaires $\ga-\gb$, pour $\ga,\gb\in S\cup\{0\}$, sont des octants de sommets dans $-\NN^\ell\times\{0_{n-\ell}\}$ et sont deux à deux emboîtés.}

\noindent[Pour plusieurs questions, on peut se contenter de la condition plus faible que les polyèdres des $\ga-\gb$, pour $\ga,\gb\in S$, sont des octants de sommets dans $-\NN^\ell\times\{0_{n-\ell}\}$, auquel cas elle ne concerne que les fibrés de rang $\geq2$. Notons aussi que l'une ou l'autre de ces propriétés est toujours satisfaite en dimension~$1$.]

On dit qu'un fibré de Higgs sur $X^o$ est \emph{sauvage et bon} le long de~$D$ si la décompo\-si\-tion~\eqref{eq:decomp} a lieu au voisinage de tout point de~$D$ après ramification locale, avec un ensemble local $\Irr(\theta)$ \emph{bon}.

De même, on dit qu'un fibré méromorphe plat $(\cV,\nabla)$ sur~$X$ à pôles le long de~$D$ admet une \emph{bonne structure formelle le long de~$D$} si, pour tout point $x\in D$, le fibré à connexion tensorisé par $\CC\lcr z_1,\dots,z_n\rcr[1/z_1\dots z_\ell]$ admet, après ramification, une décomposition \eqref{eq:LT} paramétrée par un ensemble fini et \emph{bon} $\Irr_x(\nabla)\subset\cO(*D)/\cO$. Cette propriété a été utilisée lorsque~$D$ est \emph{lisse}, dans l'étude des déformations isomonodromiques d'équations différentielles d'une variable à singularités irrégulières. Les premiers travaux dans le cas des croisements normaux, après les premiers cas considérés dans \cite{L-vdE82}, sont ceux de Majima \cite{Majima84}, poursuivis par \cite{Bibi97} dans le cas de deux variables. La situation est maintenant claire grâce à l'analyse très détaillée faite par T\ptbl Mochizuki dans \cite{Mochizuki10b,Mochizuki08}.

Comme pour les singularités régulières avec le prolongement canonique de Deligne \cite{Deligne70}, il est important de pouvoir travailler avec un $\cO_X$-module cohérent $\ccV$ tel que $\cV=\cO_X(*D)\otimes_{\cO_X}\nobreak\ccV$, appelé \emph{réseau} de~$\cV$. L'existence globale d'un tel réseau n'est pas évidente. Même dans la situation locale, lorsque $(\cV,\nabla)$ admet une bonne structure formelle, il n'est pas clair qu'existe un réseau dont le formalisé soit adapté à la décomposition \eqref{eq:LT}. L'existence d'un tel réseau est pourtant essentielle pour exhiber les propriétés asymptotiques des sections horizontales de la connexion au voisinage du diviseur: c'est d'abord dans une base locale d'un tel réseau qu'on arrive à les exprimer.

\begin{defi}[{Bon réseau \cite[\S2.3]{Mochizuki08}}]
Un \emph{bon réseau} est un $\cO_X$-sous-module cohérent sans torsion $\ccV$ du fibré méromorphe $\cV$,  qui l'engendre par tensorisation par $\cO_X(*D)$ et tel que le formalisé $\wh\ccV:=\wh\cO_{X,x}\otimes_{\cO_{X,x}}\ccV$ en chaque point $x\in D$ soit la partie invariante par l'action du groupe de Galois d'une ramification locale d'un \emph{bon réseau non ramifié}, c'est-à-dire qui se décompose de manière compatible à la décomposition \eqref{eq:LT} de $\wh\cO_{X,x}\otimes_{\cO_{X,x}}\cV$, de sorte que sur la composante $(\wh\ccV_\ga$, $\wh\nabla{}_\ga^\reg)$ soit à pôles logarithmiques au sens de \cite{Deligne70}.
\end{defi}

Malgrange \cite{Malgrange95} a montré l'existence d'un \og réseau canonique\fg qui est \emph{bon sur un ouvert de Zariski dense de~$D$}. La construction est locale et, comme pour le réseau canonique de Deligne dans le cas des singularités régulières, c'est en contrôlant les résidus le long des composantes de~$D$ que Malgrange peut globaliser diverses constructions locales. Ce réseau est appelé \emph{réseau canonique de Deligne-Malgrange} par Mochizuki. Du point de vue de l'analyse asymptotique au voisinage des singularités de $D$, ce réseau est encore insuffisant, mais Mochizuki montre:

\begin{theo}[{\cite[Cor\ptbl2.24]{Mochizuki10b}}]
Si $(\cV,\nabla)$ admet une bonne structure formelle le long de $D$, alors le réseau canonique de Deligne-Malgrange est bon en tout point de $D$.
\end{theo}

\skpt
\begin{remas}
\begin{enumeratei}
\item
La notion de bon réseau s'étend de manière évidente au cas d'une $\hb$-connexion si $\hb\neq0$. Dans le cas Higgs ($\hb=0$), il est même inutile de passer au formalisé.
\item
Dans les résultats mentionnés à la remarque \ref{rem:Higgsmodsauv}\eqref{rem:Higgsmodsauv4} ainsi qu'au théorème \ref{th:LT}, c'est la propriété \og sauvage et bon\fg qui est obtenue après éclatements, autrement dit les ensembles $\Irr_x(\theta)$ et $\Irr_x(\nabla)$ sont bons pour tout $x\in D$.
\item
Dans la suite, nous admettrons le théorème \ref{th:LT}, dans la mesure où il possède maintenant une démonstration indépendante par Kedlaya. Dans \cite{Mochizuki08}, Mochizuki ne pouvait pas se permettre ce raccourci, n'ayant alors pas à sa disposition ce théorème, et il l'a démontré par les arguments indiqués plus haut, que nous n'expliciterons pas.
\end{enumeratei}
\end{remas}

\subsection{Fibrés harmoniques sauvages (et bons)}
La propriété \og modéré\fg ou \og sauvage\fg pour un fibré plat (ou de Higgs) muni d'une métrique harmonique porte sur le champ de Higgs associé. Dans le cas d'une courbe, la propriété \og modéré\fg a été introduite par Simpson \cite{Simpson90}. Elle a été étendue par Biquard \cite{Biquard97} au cas où~$D$ est un diviseur lisse, puis par Mochizuki \cite{Mochizuki02,Mochizuki07} à celui où~$D$ est un diviseur à croisements normaux. La condition \og sauvage\fg, déjà mentionnée par Simpson \cite{Simpson90}, a été considérée sur les courbes dans \cite{Bibi98,B-B03,Bibi06b} et enfin, en toute généralité, dans \cite{Mochizuki08}.

\begin{defi}
Un fibré plat harmonique $(V,\nabla,h)$ (\resp un fibré de Higgs harmonique $(E,\theta,h)$) sur $X^o$ est \emph{sauvage} (\resp \emph{sauvage et bon}) si le fibré de Higgs associé $(E,\theta)$ est sauvage, voir \S\ref{subsec:Higgssauvage} (\resp \emph{sauvage et bon}, voir \S\ref{subsec:goodwild}).
\end{defi}

\skpt
\begin{remas}
\begin{enumeratei}
\item
Pour un fibré de Higgs \emph{harmonique}, l'holomorphie des $f_{i,k}$ (voir \S\ref{subsec:Higgssauvage}) implique la constance des $f_{i,k|D_i}$ (voir \cite[Lem\ptbl8.2]{Mochizuki07}). Dans le cas d'un fibré de Higgs harmonique satisfaisant à \eqref{eq:decomp}, il n'est pas clair \emph{a priori} que les composantes $(E_\ga,\theta_\ga-\nobreak\rd\ga\otimes\nobreak\id)$, munies de la métrique induite, soient harmoniques. Il faut donc les traiter comme des fibrés de Higgs quelconques, et imposer la constance des $f_{i,k|D_i}$.

\item
Mochizuki \cite[Chap\ptbl8]{Mochizuki07} donne, pour un fibré de Higgs harmonique, un critère de modération par restriction aux courbes transverses à la partie lisse de~$D$, réminiscent de celui donné par Deligne \cite{Deligne70} pour les connexions plates à singularités régulières.
\end{enumeratei}
\end{remas}

\section{Correspondance de Hitchin-Kobayashi sauvage}\label{sec:HK}

Dans la suite du texte, nous allons esquisser la démonstration du fait que tout $\cD_X$\nobreakdash-module holonome simple provient d'un $\cD$-module  holonome avec structure de twisteur polarisable (point \eqref{enum:programme3} du \S\ref{sec:strategie}). Nous commençons par la fin, à savoir le point \eqref{enum:programme6}.

Dans la situation globale (voir convention \ref{subsec:convention}), soit $(V,\nabla)$ un fibré algébrique plat sur~$X^o$. Il correspond de manière biunivoque à un $\cO_X(*D)$-module cohérent~$\cV$ à connexion plate~$\nabla$. Supposons $(V,\nabla)$ simple. Mochizuki montre l'existence d'une métrique harmonique pour $(V,\nabla)$ avec de bonnes propriétés. Décrivons les étapes en renvoyant plus bas pour les définitions précises.

\begin{enumeratea}
\item
La construction de Malgrange d'un réseau canonique \cite{Malgrange95} permet de munir $(\cV,\nabla)$ d'une filtration parabolique $_\bbullet\cV^\DM$.
\item
Étant donné un fibré ample $L$ sur~$X$, on associe à tout $\cO_X(*D)$-module cohérent plat parabolique une pente $\mu_L$, d'où une notion de $\mu_L$-stabilité. Alors $(\cV,\nabla)$ est simple si et seulement si $(\cV,{}_\bbullet\cV^\DM,\nabla)$ est $\mu_L$-stable.
\item
On a aussi une notion de nombres caractéristiques paraboliques. D'après le théorème \ref{th:LT} (en utilisant la version de Kedlaya), quitte à changer~$X$ et~$D$, on peut même supposer que $(\cV,\nabla)$ admet une bonne structure formelle le long de~$D$. Dans ce cas, les nombres caractéristiques de $(\cV,{}_\bbullet\cV^\DM)$ sont nuls.
\item
La correspondance de Hitchin-Kobayashi consiste alors, dans ce cadre, en la construction d'une métrique harmonique~$h$ pour $(V,\nabla)$ adaptée à la filtration $({}_\bbullet\cV^\DM)$, d'où un fibré de Higgs harmonique $(E,\theta,h)$ et une variation de structure de twisteur polarisée pure de poids~$0$ sur $X^o$. Cette construction repose aussi sur le fait que le fibré de Higgs harmonique $(E,\theta,h)$ ainsi obtenu est sauvage et bon.
\end{enumeratea}

Les questions qui restent en suspens sont de savoir si $(E,\theta)$ se prolonge (et comment) à~$X$, et si la variation de structure de twisteur polarisée se prolonge (et comment) à~$X$. Nous les aborderons au \S\ref{sec:prolong}.

\subsection{Filtrations paraboliques et métriques adaptées}\label{subsec:locab}
Nous nous plaçons dans la situation locale ou globale (convention \ref{subsec:convention}). Mochizuki a considéré une propriété analogue à la propriété suivante dans \cite[\S4.2]{Mochizuki07}.\enlargethispage{\baselineskip}%

\subsubsection*{Filtration parabolique}
Soit $_\inftyg E$ un $\cO_X(*D)$-module cohérent sans torsion. Une \emph{filtration parabolique} de $_\inftyg E$ consiste en la donnée d'une filtration croissante, indexée par~$\RR^\ccI$ muni de son ordre partiel naturel, de $_\inftyg E$ par des sous-$\cO_X$-modules cohérents sans torsion $_\bma E$, qui satisfait aux propriétés suivantes:
\begin{enumerate}
\item (translation)
pour tout $\bma\in\RR^\ccI$, on a $_\inftyg E=\cO_X(*D)\otimes_{\cO_X}{}_\bma E$ et, pour tout $\bnn\in\ZZ^\ccI$, $_{\bma-\bnn}E=\cO_X(-\sum_{i\in\ccI}n_iD_i)\otimes_{\cO_X}{}_\bma E$;
\item (finitude)
il existe pour tout $i\in\ccI$ un sous-ensemble $\ccA_i\subset\RR$ fini modulo $\ZZ$ tel que la filtration soit déterminée par sa restriction à $\ccA=\prod_{i\in\ccI}\ccA_i$, c'est-à-dire que pour tout $\bma'\in\RR^\ccI$, on a $_{\bma'}E=\bigcup_{\substack{\bma\in\ccA\\\bma\leq\bma'}}{}_\bma E$.
\end{enumerate}

Si les $_\bma E$ sont $\cO_X$-localement libres, on dira aussi que $({}_\inftyg E,{}_\bbullet E)$ est un \emph{fibré méromorphe parabolique sur $(X,D)$}. Si ${}_\inftyg E$ est localement libre de rang $1$, la donnée d'une filtration parabolique est équivalente à la donnée de $\bmb\in\RR^\ccI$ modulo $\ZZ^\ccI$. On a alors localement $_\bma E\simeq\cO_X(\sum_i[a_i+b_i]D_i)$. La proposition suivante simplifie différentes notions introduites dans \cite[Chap\ptbl3]{Mochizuki06}, \cite[Chap\ptbl4]{Mochizuki07} et \cite[\S2]{I-S07}.

\begin{prop}[{\cite[Th\ptbl2.4.20]{Borne09}, \cite[Th\ptbl4.2]{H-S08}}]\label{prop:locab}
Tout fibré méromorphe parabolique $({}_\inftyg E,{}_\bbullet E)$ sur $(X,D)$ est localement abélien, \ie localement isomorphe à une somme directe de fibrés paraboliques de rang $1$.
\end{prop}

\begin{defi}\label{def:eprime}
Soit $({}_\inftyg E,{}_\bbullet E)$ un fibré méromorphe parabolique sur $(X,D)$. Une base locale $\bme$ de $_\inftyg E$ comme $\cO_X(*D)$-module est dite \emph{adaptée} à la filtration parabolique si elle définit une décomposition de $({}_\inftyg E,{}_\bbullet E)$ en fibrés paraboliques de rang $1$. Chaque élément $e_k$ a alors un multi-ordre $\bma(k)$. On associe à $\bme$ une base normalisée $\bme'$ définie par $e'_k=\prod_{i=1}^\ell|z_i|^{a_i(k)}\cdot e_k$.
\end{defi}

\subsubsection*{Nombres caractéristiques paraboliques}
Soit $({}_\inftyg E,{}_\bbullet E)$ un fibré méromorphe parabolique.  Considérons l'un des fibrés $_\bma E$ (pour $\bma\in\ccA$, c'est suffisant). Sur chaque composante $D_i$ de~$D$ on dispose alors du $\cO_{D_i}$-module localement libre $_\bma E/{}_{\bma^{-_i}} E$, si $\bma^{-_i}$ est le prédécesseur de $\bma$ dans la direction $i$ uniquement. Notons $_\bma E_{||D_i}$ ce fibré et $\rg{}_\bma E_{||D_i}$ son rang. On observe que, pour $\bmb\in\ccA$, la restriction usuelle $_\bmb E_{|D_i}={}_\bmb E/{}_{\bmb-1_i}E$ a pour rang
\[
\rg{}_\bmb E_{|D_i}=\sum_{\substack{a_i\in{}]b_i-1,b_i]\\ a_j=b_j\,\forall j\neq i}}\rg{}_\bma E_{||D_i}
\]

\begin{defi}[{\cite[\S3.1.2]{Mochizuki06}}]
Dans la situation globale, soit $L$ un fibré ample sur~$X$. Le \emph{degré parabolique} $\parc\deg_L({}_\inftyg E,{}_\bbullet E)$ est défini par la formule, indépendante du choix de $\bmb\in\RR^\ccI$,
\[
\parc\deg_L({}_\inftyg E,{}_\bbullet E)=\deg_L{}_\bmb E-\sum_{i\in\ccI}\Big(\sum_{a_i\in{}]b_i-1,b_i]}a_i\rg{}_\bma E_{||D_i}\Big)\deg_LD_i.
\]
La \emph{pente} $\mu_L({}_\inftyg E,{}_\bbullet E)$ est le quotient $\parc\deg_L({}_\inftyg E,{}_\bbullet E)/\rg{}_\inftyg E$.
\end{defi}

On peut aussi définir une classe $\parc c_1({}_\inftyg E,{}_\bbullet E)$ en remplaçant dans la formule ci-dessus le nombre $\deg_LD_i$ par la classe $[D_i]$ dans $H^2(X,\RR)$ et $\deg_L{}_\bmb E$ par $c_1({}_\bmb E)$, et on peut aussi définir un nombre $\parc\deg_L\mathrm{ch}_2({}_\inftyg E,{}_\bbullet E)$ (voir \loccit).

\subsubsection*{Métrique hermitienne adaptée à une filtration parabolique}
Dans la situation locale, soit $({}_\inftyg E,{}_\bbullet E)$ un fibré méromorphe parabolique sur $(X,D)$. Soit par ailleurs~$h$ une métrique hermitienne sur $E={}_\inftyg E_{|X^o}$. Pour chaque $\bma\in\RR^\ccI$, définissons le sous-faisceau de $\cO_X$-modules $_\bma \wt E\subset j_*E$ par
\[
\forall U\subset X,\quad{}_\bma \wt E(U)=\Big\{e\in E(U\moins D)\mid\forall\epsilon>0,\,|e|_h=O\big(\ts\prod_{i\in\ccI}|z_i|^{-a_i-\epsilon}\big)\,\text{loc.\,sur\,$U$}\Big\},
\]
et $_\inftyg\wt E=\bigcup_\bma {}_\bma \wt E$, qui est un $\cO_X(*D)$-module, filtré par les  $\cO_X$-sous-modules sans torsion~$_\bma\wt E$. En général, ces faisceaux n'ont aucune propriété de cohérence.

\begin{defi}[{\cite[\S3.5]{Mochizuki06}}]\label{def:adapte}
La métrique~$h$ est dite \emph{adaptée au fibré méro\-morphe parabolique $({}_\inftyg E,{}_\bbullet E)$} si ${}_\bma\wt E={}_\bma E$ pour tout $\bma\in\RR^\ccI$.
\end{defi}

\subsection{La filtration de Deligne-Malgrange}
Revenons à notre problème, dans la situation globale. Soit $(\cV,\nabla)$ un $\cO_X(*D)$-module cohérent à connexion intégrable. Lorsque $(\cV,\nabla)$ est à \emph{singularités régulières} le long de~$D$, Deligne \cite{Deligne70} a construit un fibré vectoriel canonique, sur lequel la connexion est à pôles logarithmiques et les valeurs propres $\alpha_i$ de l'endomorphisme résidu sur la composante $D_i$ de~$D$ ont une partie réelle dans $[0,1[$. Pour chaque $\bma\in\RR^\ccI$, on peut définir le fibré à connexion logarithmique~$_\bma\cV$ en imposant que $-\reel\alpha_i\in{}]a_i-1,a_i]$ pour tout $i\in\ccI$. On obtient ainsi la \emph{filtration canonique de Deligne} et un fibré méromorphe plat parabolique $(\cV,{}_\bbullet\cV,\nabla)$, tout ceci compatible à la formation du déterminant. Dans le cas de rang $1$, la connexion sur le fibré $C^\infty$ associé est à coefficients distributions, et la formule de Chern-Weil (au sens des courants) pour $c_1({}_\bmb\cV)$  montre que $\parc c_1(\cV,{}_\bbullet\cV)=\nobreak0$, et ceci reste vrai en tout rang par passage au déterminant, de même que l'égalité $\parc\deg_L(\cV,{}_\bbullet\cV)=0$.

On remarque aussi que tout $\cO_X(*D)$-sous-module cohérent de $\cV$ stable par la connexion est encore $\cO_X(*D)$-localement libre, et la connexion y est encore à singularités régulières. De plus, la filtration de Deligne de $\cV$ induit sur ce sous-module sa propre filtration de Deligne.

On dira que le fibré méromorphe plat parabolique $(\cV,{}_\bbullet\cV,\nabla)$ est $\mu_L$-stable si tout sous-fibré méromorphe plat, muni de la filtration parabolique induite, est de pente strictement plus petite.

On vérifie alors que $(\cV,{}_\bbullet\cV,\nabla)$ est $\mu_L$-stable si et seulement si $(\cV,\nabla)$ est simple (en effet, la filtration de Deligne induit sur tout sous-fibré méromorphe plat non trivial la filtration de Deligne de celui-ci, qui est donc de pente nulle, en contradiction avec la stabilité).

Lorsque $(\cV,\nabla)$ est à singularités irrégulières, Malgrange \cite{Malgrange95} a construit un réseau canonique (\cf \S\ref{subsec:goodwild}) en imposant des conditions analogues sur les résidus des connexions~$\wh\nabla_\ga^\reg$, qui existent \emph{a priori} sur un ouvert de Zariski dense de chaque $D_i$. On en déduit une filtration canonique $({}_\bbullet\cV^\DM,\nabla)$, dite \emph{de Deligne-Malgrange}. Il faut noter que chaque $({}_\bma\cV^\DM,\nabla)$ n'est pas nécessairement logarithmique, et n'est pas nécessairement localement libre. Néanmoins, chaque ${}_\bma\cV^\DM$ est $\cO_X$-cohérent et réflexif (voir \cite[Lem\ptbl2.7.8]{Mochizuki08}).

Si de plus $(\cV,\nabla)$ admet une bonne structure formelle le long de~$D$ alors, comme indiqué au \S\ref{subsec:goodwild}, le réseau canonique de Deligne-Malgrange est localement la partie invariante dans une ramification d'un réseau qui se décompose formellement au voisinage de chaque point de~$D$ comme $(\cV,\nabla)$. Il en résulte que c'est un fibré vectoriel, et $(\cV,{}_\bbullet\cV^\DM)$ est un fibré méromorphe plat parabolique. Dans le cas de rang $1$, se préoccuper de bonté et de ramification est superflu et, utilisant le quasi-isomorphisme $\Omega^\cbbullet(\log D)\simeq\Omega^\cbbullet(*D)$ (voir \cite[Prop\ptbl II.3.13]{Deligne70}), on montre que la connexion $\nabla+\ov\partial$ sur $\cC^\infty_X\otimes\cV$ s'écrit comme la somme d'une connexion plate $C^\infty$ logarithmique et d'une forme exacte $\rd\varphi$ avec $\varphi\in\Gamma(X,\cC^\infty_X(*D))$. Il s'ensuit que, comme dans le cas logarithmique, on a, en tout rang, $\parc c_1(\cV,{}_\bbullet\cV)=\nobreak0$ et $\parc\deg_L(\cV,{}_\bbullet\cV)=0$.

On montre de même que tout $\cO_X(*D)$-sous-module cohérent de $\cV$ stable par la connexion est encore $\cO_X(*D)$-localement libre, et admet une bonne structure formelle le long de~$D$ (avec des facteurs exponentiels locaux contenus dans ceux de $(\cV,\nabla)$). Enfin, la filtration de Deligne-Malgrange de $\cV$ induit celle de ses sous-modules.

On en déduit, de même que plus haut, l'équivalence $\mu_L$\nobreakdash-stabilité de $(\cV,{}_\bbullet\cV^\DM,\nabla)\iff{}$ simplicité de $(\cV,\nabla)$ (voir \cite[\S2.7.2.2]{Mochizuki08}). On a aussi:

\begin{prop}[{\cite[Cor\ptbl14.3.4]{Mochizuki08}}]
On a $\parc\deg_L\mathrm{ch}_2({}_\bbullet\cV^\DM)=0$.
\end{prop}

\subsection{Construction d'une métrique harmonique adaptée}
On suppose que $(\cV,\nabla)$ est simple et admet une bonne structure formelle le long de~$D$. On cherche à construire une métrique harmonique adaptée à la structure parabolique de Deligne-Malgrange.

\subsubsection*{Le cas de rang $1$}
Il est instructif de commencer par considérer le cas des fibrés de rang~$1$. On remarque d'abord que, en rang $1$, la pseudo-courbure $G(\nabla,h)$ (\cf\eqref{eq:pseudocourbure}) est égale à deux fois la courbure de~$h$ (voir \cite[Lem\ptbl2.31]{Mochizuki09c}). Il s'agit donc de construire une métrique à courbure nulle adaptée à une filtration parabolique, elle-même déterminée par la donnée de $\bma\in\RR^\ccI$. On construit une métrique singulière~$h_0$ sur~${}_\bma\cV$ en imposant d'abord que, dans toute situation locale, une base locale $e$ de~$_\bma\cV$ ait pour norme $\Vert e\Vert_{|h_0}=|z|^{-\bma}\cdot\Vert e\Vert_{h_\loc}$, où~$h_{\loc}$ est une métrique $C^\infty$ locale sur ${}_\bma\cV$, puis en utilisant une partition de l'unité pour recoller. La courbure $R(h_0)$ est la somme d'une forme $C^\infty$ fermée $R'(h_0)$ de type $(1,1)$ sur~$X$ et d'un courant fermé de type $(1,1)$ porté par~$D$, et sa classe est $c_1({}_\bma\cV)$. La correction parabolique est faite pour que la classe de $R'(h_0)$ soit égale à $\parc c_1(\cV,{}_\bbullet\cV)$, dont on a vu plus haut qu'elle est nulle. La théorie de Hodge implique alors que $R'(h_0)=\ov\partial\partial g$ pour une certaine fonction $g$ de classe $C^\infty$, et la métrique $h=e^{-g}h_0$ est encore adaptée à ${}_\bma\cV$ et de courbure nulle.

\subsubsection*{Le cas des courbes}
Dans la situation globale, soit~$X$ une courbe algébrique lisse et $(\cV,\nabla)$ un fibré méromorphe à connexion sur $(X,D)$. La filtration de Deligne-Malgrange est définie de manière classique dans ce cas (voir \cite{Malgrange95}) et la condition de bonté est trivialement satisfaite. Supposons $(\cV,\nabla)$ simple. Il résulte essentiellement des travaux de Simpson \cite{Simpson90} (voir \cite{Bibi98}) qu'il existe une métrique harmonique~$h$ pour $(\cV_{|X^o},\nabla)$ adaptée à la filtration de Deligne-Malgrange.

Dans \cite{B-B03} Biquard et Boalch étendent ce résultat à une situation parabolique plus générale, et le précisent en une correspondance de Hitchin-Kobayashi entre fibrés plats et fibrés de Higgs.

Dans \cite[\S13.4]{Mochizuki08}, Mochizuki donne une autre démonstration de ce résultat, dans l'esprit de celle de \cite{Simpson90}, et il précise une propriété d'unicité de la métrique harmonique.

\subsubsection*{Un théorème de type Metha-Ramanathan}
Ce résultat de réduction aux courbes générales lorsque $\dim X\geq2$ est important à plusieurs endroits de la preuve qui suit. L'argument de Simpson \cite{Simpson92} a déjà été généralisé par Mochizuki dans \cite{Mochizuki06} au cas des singularités régulières, et la démonstration est étendue au cas des singularités irrégulières:

\begin{prop}[{\cite[Cor\ptbl13.2.3]{Mochizuki08}}]\label{prop:MR}
Dans la situation globale, soit $(\cV,\nabla)$ un fibré méromorphe plat à pôles le long de~$D$ et $L$ un fibré ample sur $X$. Alors $(\cV,\nabla)$ est simple si et seulement si sa restriction à toute courbe intersection complète assez générale de sections de~$L^{\otimes m_\nu}$, pour une suite $m_\nu\to+\infty$, est simple.
\end{prop}

\subsubsection*{Construction d'une métrique harmonique adaptée }
Soit $(\cV,\nabla)$ un fibré méromorphe plat sur $(X,D)$. Supposons que $(\cV,\nabla)$ admette une \emph{bonne structure formelle le long de~$D$}, et donne ainsi lieu à un fibré méromorphe plat parabolique $(\cV,{}_\bbullet\cV^\DM,\nabla)$. Supposons aussi que $(\cV,\nabla)$ est \emph{simple}, de sorte que $(\cV,{}_\bbullet\cV^\DM,\nabla)$ est $\mu_L$-stable, à nombres caractéristiques nuls. De plus, comme indiqué plus haut, la première classe de Chern parabolique est nulle. On peut définir un fibré méromorphe plat parabolique de rang~$1$, à savoir le déterminant $\det(\cV,{}_\bbullet\cV^\DM,\nabla)$, qui est aussi à $\parc c_1$ nul et qui, d'après ce qu'on a vu plus haut, admet une métrique harmonique~$h_{\det}$.

\begin{theo}[{\cite[Th\ptbl25.28]{Mochizuki07}, \cite[Th\ptbl5.16]{Mochizuki09c}, \cite[Th\ptbl16.1.1]{Mochizuki08}}]\label{th:existharm}
Dans ces conditions, il existe une unique métrique harmonique adaptée à $(\cV,{}_\bbullet\cV^\DM,\nabla)$ normalisée par $\det(h)=h_{\det}$. De plus cette métrique fait du fibré de Higgs harmonique associé $(E,\theta,h)$ un fibré de Higgs sauvage et bon.
\end{theo}

Considérons d'abord la dernière propriété, qui sera utile pour montrer l'unicité, et en dernier ressort l'existence en dimension $\geq3$. C'est un résultat plus général.

\begin{prop}[l'adaptation implique la sauvage bonté, {\cite[Prop\ptbl13.5.2]{Mochizuki08}}]
Soit $(\cV,\nabla)$ un fibré méromorphe plat sur $(X,D)$ admettant une bonne structure formelle le long de~$D$. Supposons aussi $(\cV,\nabla)_{|X^o}$ muni d'une métrique harmonique~$h$ adaptée à ${}_\bbullet\cV^\DM$. Alors le fibré de Higgs harmonique associé $(E,\theta,h)$ est sauvage et bon.
\end{prop}

\subsubsection*{Unicité dans \ref{th:existharm}}
Prenons deux telles métriques $h_1$ et $h_2$. En restriction à une courbe générale $C$ comme en \ref{prop:MR}, $(\cV,\nabla)$ est simple et sauvage (automatiquement bon en dimension $1$). Anticipant sur la section suivante (Théorèmes \ref{th:acceptable} et \ref{th:harmacceptable}), l'estimation du théorème \ref{th:acceptable}\eqref{th:acceptable2} pour $({}_\bbullet\cV^\DM,h)$ se restreint à la courbe $C$, ce qui permet de voir que les restrictions de $h_1,h_2$ à $C^o$ sont adaptées à $({}_\bbullet\cV^\DM)_{|C}$. L'unicité vue plus haut dans le cas des courbes montre que $h_{1|C^o}=h_{2|C^o}$. Puisqu'on peut faire passer une telle courbe générale par chaque point de $X^o$, on en déduit l'unicité.

\subsubsection*{Existence dans \ref{th:existharm}}
La technique a été développée par Mochizuki dans le cas modéré dans \cite{Mochizuki07,Mochizuki09c} et s'adapte au cas sauvage et bon à l'aide des résultats déjà obtenus et ceux des paragraphes \ref{subsec:accept} et \ref{subsec:acceptwild} ci-dessous. Résumons-la \emph{très} rapidement.

D'une part Mochizuki étend des résultats de Donaldson et Simpson (voir \cite{Simpson88}), et d'autre part il commence par le cas des surfaces. Dans ce cas, il construit une famille, paramétrée par $\epsilon>0$, de perturbations de la structure parabolique pour supprimer l'éventuelle partie nilpotente des gradués paraboliques des résidus de $\wh\nabla_\ga^\reg$ le long des composantes de~$D$. À l'aide des résultats de Donaldson et Simpson généralisés, et à partir d'une métrique convenable sur les différents fibrés gradués paraboliques plats sur les composantes $D_i$, il obtient pour chaque $\epsilon$ une métrique harmonique\footnote{au premier sens de la note page \pageref{footnote:harmonique}.} adaptée à la filtration parabolique perturbée. Il montre ensuite la convergence de ces métriques pour $\epsilon\to0$, en un sens convenable, vers une métrique harmonique\footnote{au sens \og pluri-harmonique\fg de la note page \pageref{footnote:harmonique}.} adaptée à la filtration parabolique $({}_\bbullet\cV^\DM)$.

Le cas de la dimension $\geq3$ se traite grâce à l'argument d'unicité: la métrique définie sur chaque surface assez générale, grâce à \ref{prop:MR}, est bien la restriction d'une métrique existant sur $X^o$, et elle satisfait aux propriétés requises.

\section{Prolongement de fibrés harmoniques sauvages}\label{sec:prolong}

\subsection{Les problèmes du prolongement}
On se place dans la situation locale. Soit $(E,\theta,h)$ un fibré de Higgs harmonique sauvage et bon sur $X^o$. On cherche en particulier à résoudre les problèmes suivants:
\begin{enumeratea}
\item\label{enum:pbextension1}
prolonger le fibré $E$ en un fibré méromorphe $\cE$ (\ie $\cO_X(*D)$-module localement libre de rang fini) sur~$X$,
\item\label{enum:pbextension2}
montrer que $\theta$ se prolonge de manière méromorphe au prolongement $\cE$ (\ie les coefficients de $\theta$ dans une base locale de $\cE$ sont méromorphes),
\item\label{enum:pbextension3}
pour chaque $\hb$, prolonger aussi les fibrés plats $(V^\hb,\nabla)$ en des fibrés méromorphes plats $(\cV^\hb,\nabla)$.
\item\label{enum:pbextension4}
effectuer ce dernier prolongement de manière holomorphe par rapport à $\hb$.
\end{enumeratea}
On explique dans cette section comment Mochizuki procède dans le cas sauvage et bon \cite{Mochizuki08}, après avoir résolu ces questions dans le cas modéré \cite{Mochizuki07}.

Dans la mesure où le problème de prolongement local est résolu canoniquement à l'aide de la métrique (voir ci-dessous), il conduit à des résultats applicables dans la situation globale.

\subsection{Fibrés holomorphes hermitiens acceptables}\label{subsec:accept}
Considérons la situation locale. Soit $(E,h)$ un fibré holomorphe hermitien sur $X^o$. La condition d'\emph{acceptabilité} remonte à l'article de Cornalba et Griffiths \cite{C-G75}, et a été utilisée dans ce cadre par Simpson \cite{Simpson88,Simpson90}:

\begin{defi}
Le fibré hermitien $(E,h)$ est \emph{acceptable} si la norme de la courbure de~$h$, calculée par rapport à~$h$ elle-même (sur le fibré des endomorphismes) et à la métrique de Poincaré sur $U\moins D$, est \emph{bornée}.
\end{defi}

Mochizuki raffine les résultats sur les fibrés acceptables de la manière suivante, avec la notation utilisée à la définition \ref{def:adapte}:

\skpt
\begin{theo}[{\cite{Mochizuki07} et \cite[21.3.1--3]{Mochizuki08}}]\label{th:acceptable}
\begin{enumeratei}
\item\label{th:acceptable1}
Si $(E,h)$ est acceptable, alors chaque $_\bma\wt E$ est $\cO_X$-localement libre et $({}_\bma\wt E)$ fait de~$_\inftyg\wt E$ un fibré méromorphe parabolique $(_\inftyg\wt E,{}_\bbullet\wt E)$.
\item\label{th:acceptable2}
De plus, si $\bme$ est une base locale de $_\inftyg\wt E$ adaptée à la décomposition donnée par la proposition \ref{prop:locab} et $\bme'$ la base normalisée associée (voir définition \ref{def:eprime}), les valeurs propres~$\eta(z)$ de la matrice $h(\bme',\bme')$ de la métrique dans cette base satisfont aux inégalités
\[
C\Big(\sum_{i=1}^\ell\rL(z_i)\Big)^{-N}\leq\eta(z)\leq C'\Big(\sum_{i=1}^\ell\rL(z_i)\Big)^N
\]
pour $C,C',N$ positifs convenables, en posant $\rL(z)=|{\log|z|}|$ ($|z|<1$).
\item\label{th:acceptable3}
Enfin, le fibré $(\End{}_\inftyg\wt E,h)$ est aussi acceptable, et $_{\mathbf{0}\!}\End{}_\inftyg\wt E$ est le faisceau des endomorphismes de $_\inftyg\wt E$ qui préservent la filtration $_\bma\wt E$ et dont la restriction à chaque composante $D_i$ préserve la filtration naturelle de $_\bma\wt E_{|D_i}$.
\end{enumeratei}
\end{theo}

\subsection{Acceptabilité des fibrés harmoniques sauvages et bons}\label{subsec:acceptwild}
Restons dans un cadre local comme plus haut. Soit $(E,\theta,h)$ un fibré de Higgs harmonique sur $X^o$. La courbure de la connexion de Chern se calcule, du fait de l'harmonicité, par la formule $R(h,\ov\partial_E)=-[\theta,\theta^\dag]$. Plus généralement, pour tout $\hb\in\CC$, $R(h,\ov\partial_E+\hb\theta^\dag)=-(1+|\hb|^2)[\theta,\theta^\dag]$. Lorsque $(E,\theta,h)$ est sauvage et bon, Mochizuki montre l'acceptabilité de $(E,h)$, qui entraîne donc celle de tous les $(V^\hb,h)$. Plus précisément, il  donne une interprétation géométrique de cette propriété, que nous décrivons maintenant, généralisant celle donnée par Simpson \cite{Simpson90} en dimension $1$ et dans la situation modérée.

La décomposition \eqref{eq:decomp} peut se raffiner: $(E,\theta)\simeq\bigoplus_{(\ga,\alphag)}(E_{(\ga,\alphag)},\theta_{(\ga,\alphag)})$, en utilisant \eqref{eq:decompvp}. Cette décomposition ne dépend que du champ de Higgs, et on va l'analyser par rapport à la métrique~$h$. On dira que $E_{(\ga,\alphag)}$ est \emph{$g_{(\ga,\alphag),(\gb,\betag)}$-asymptotiquement $h$-orthogonal à $E_{(\gb,\betag)}$} si la norme de la projection de $E_{(\gb,\betag)}$ sur $E_{(\ga,\alphag)}$ parallèlement à $E_{(\ga,\alphag)}^{\perp_h}$ et celle de $E_{(\ga,\alphag)}$ sur $E_{(\gb,\betag)}$ parallèlement à $E_{(\gb,\betag)}^{\perp_h}$ sont localement bornées par une fonction $g_{(\ga,\alphag),(\gb,\betag)}(z)$.

Dans la suite, on prendra $g_{\epsilon,(\ga,\alphag),(\gb,\betag)}(z)=\exp(-\epsilon|z^{\ord(\ga-\gb)}|)\prod_{j|\alpha_i\neq\beta_i}|z_i|^\epsilon$ ($\epsilon>0$). La condition (Bon) implique en effet que pour tous $\ga\neq\gb$, $\ga-\gb=z^{-\bmm}\gc(z)$ avec~$\gc$ holo\-mor\-phe et $\gc(0)\neq0$, pour un certain multi-indice $\bmm\in\ZZ^\ell\moins\NN^\ell$, qu'on note $\ord(\ga-\nobreak\gb)$. On voit donc que si $\ga\neq\gb$, la fonction $g_{\epsilon,(\ga,\alphag),(\gb,\betag)}$ est à décroissance exponentielle au voisinage de l'origine, tandis que si $\ga=\gb$, on a seulement une décroissance modérée.

\begin{theo}[{\cite[Chap\ptbl8]{Mochizuki07} et \cite[Chap\ptbl7]{Mochizuki08}}]\label{th:harmacceptable}
Soit $(E,\theta,h)$ un fibré de Higgs harmonique sauvage et bon. Alors $(E,h)$ est acceptable, de même que tout $(V^\hb,h)$, et $(\End V^\hb,h)$. Plus précisément, pour tous $(\ga,\alphag),(\gb,\betag)$, 
\begin{enumeratei}
\item
$E_{(\ga,\alphag)}$ est $g_{\epsilon,(\ga,\alphag),(\gb,\betag)}$-asymptotiquement $h$-orthogonal à $E_{(\gb,\betag)}$ pour $\epsilon>0$ assez petit;
\item
la norme de la composante $[\theta,\theta^\dag]_{(\ga,\alphag),(\gb,\betag)}$ relativement à~$h$ et à la métrique de Poincaré est $O(g_{\epsilon,(\ga,\alphag),(\gb,\betag)})$.
\end{enumeratei}
\end{theo}

\begin{rema}
Ce type d'estimation remonte au travaux de Simpson \cite[\S2]{Simpson90} et est pour cette raison baptisé \og Simpson's main estimate\fg par Mochizuki. Il est remarquable que la présence de parties polaires non nulles \emph{améliore nettement} les estimations d'orthogonalité asymptotique par rapport au cas modéré, puisqu'exponentiellement petites. Un tel phénomène avait déjà été observé en dimension $1$ par Biquard et Boalch \cite[Lem\ptbl4.6]{B-B03}. Néanmoins, ne nous réjouissons pas trop vite...
\end{rema}

\subsection{Prolongement à $\hb$ fixé}
Supposons toujours $(E,\theta,h)$ harmonique sauvage et bon dans la situation locale du \S\ref{subsec:Higgssauvage}. Les résultats généraux sur les fibrés hermitiens acceptables (théorème \ref{th:acceptable}) et le théorème d'acceptabilité des fibrés harmoniques sauvages et bons (théorème~\ref{th:harmacceptable}) permettent de répondre aux questions \eqref{enum:pbextension1} et \eqref{enum:pbextension2} du début de cette section. Plus précisément, pour chaque $\hb$, on obtient un fibré méromorphe parabolique sur~$X$, noté $(\ccP\cE^\hb,\ccP{}_\bbullet \cE^\hb)$.

\begin{theo}[{\cite[Cor\ptbl8.89]{Mochizuki07}, \cite[Th\ptbl7.4.5]{Mochizuki08}}]\label{th:nablahbmero}
Pour tout $\hb$, la $\hb$\nobreakdash-con\-ne\-xion $\nabla^\hb$ est méromorphe sur $\ccP{}_\bma \cE^\hb$ pour tout $\bma\in\RR^\ell$ (et logarithmique dans le cas modéré), et en fait un bon réseau, qui est non ramifié si $\theta$ l'est.
\end{theo}

La démonstration du cas modéré peut être adaptée et étendue au cas sauvage, grâce notamment à l'orthogonalité asymptotique vue ci-dessus. Supposons de plus que $(E,\theta,h)$ est sauvage sans ramification et bon. La décomposition \eqref{eq:decomp} fait intervenir un ensemble fini $\Irr(\theta)$ de parties polaires. Mochizuki obtient de plus:
\par\smallskip
\noindent(\ref{th:nablahbmero}\,$*$)\enspace
\emph{Pour tout $\hb\neq0$, l'ensemble $\Irr(\nabla^\hb)$ paramétrant la décomposition \eqref{eq:LT} pour $(\ccP\cE^\hb,\nabla^\hb)$ est égal à $(1+|\hb|^2)\Irr(\theta)$.}
\par\smallskip
Ce comportement complète le comportement des valeurs propres du résidu de $(\wh\nabla^\hb)_{(1+|\hb|^2)\ga}^\reg$, qui est régi par la fonction $\ge$ (\og eigenvalue\fg) introduite dans \eqref{eq:pe}, ainsi que l'avaient montré Simpson en dimension $1$ et Mochizuki en toute dimension, dans le cas modéré. La fonction $\gp$ (\og parabolique\fg) régit, quant à elle, le comportement par rapport à $\hb$ de la structure parabolique, que nous n'avons pas détaillé ici. La mauvaise nouvelle, anticipée dans l'exemple \ref{exem:rankone}, est que, contrairement à $\ge$, le comportement de $\Irr(\nabla^\hb)$ n'est pas holomorphe en~$\hb$.

\subsection{Prolongement à $\hb$ variable}\label{subsec:hbvariable}
Dans le cas modéré, le passage de $\hb$ fixé à $\hb$ variable n'engendre pas de complication importante. Il faut simplement prendre garde à l'uniformité locale par rapport à $\hb$ des estimations asymptotiques. En particulier, dans l'estimation du théorème \ref{th:acceptable}\eqref{th:acceptable2}, on remplace les $\big(\sum\rL(z_i)\big)^{\pm N}$ par $\prod|z_i|^{\mp\epsilon}$. Aussi je n'insisterai pas sur ce point, traité en détail dans \cite{Mochizuki07}.

Par contre, dans le cas sauvage, le comportement de l'exemple \ref{exem:rankone} vu à la fin du~\S\ref{sec:dictionnaire} n'est pas du tout anodin, et nécessite la mise en place d'une description fidèle, dans le cas sauvage et bon à plusieurs variables, des fibrés méromorphes munis d'une $\hb$-connexion intégrable en terme de l'objet formel associé et d'une \emph{structure de Stokes}. C'est ce qui est fait dans les chapitres 2 à 4, et 20 (appendice) de \cite{Mochizuki08}, et représente un apport majeur dans la théorie des $\cD$-modules holonomes, indépendamment des autres points considérés ici. La question qui nous importe est alors traitée dans les chapitres 9 à 11. Le tout occupant plus de 200 pages, il ne sera pas question de rentrer dans les détails.

Considérons la situation simplifiée, à $\hb=1$ fixé, d'un fibré méromorphe $(\cV,\nabla)$ à connexion sur un disque $\Delta$ de coordonnée $z$, avec un unique pôle en $z=0$. Le théorème de Levelt-Turrittin donne, après une ramification convenable $z'\mto z=z^{\prime q}$, une décomposition du formalisé de $(\cV,\nabla)$ en~$0$. Supposons pour simplifier l'explication que $q=1$ (cas non ramifié). On a donc une décomposition \eqref{eq:LT}.

\begin{prop}\label{prop:isomono}
Soit $t$ un réel ${}>0$. Il existe, de manière canonique et fonctorielle, un fibré méromorphe à connexion $(\cV_t,\nabla_t)$ qui a pour décomposition formelle
\[
(\wh\cV_t,\wh\nabla_t)\simeq\bigoplus_{\ga}(\wh\cV_\ga,\wh\nabla{}_\ga^\reg+t\cdot\rd\ga\otimes\id_{\wh\cV_\ga}).
\]
\end{prop}

En particulier, pour $t=1$ on retrouve $(\cV,\nabla)$, et si $t=|\tau|$ avec $\tau\in\CC^*$, alors la famille $(\cV_t,\nabla_t)$ est \emph{isomonodromique} par rapport à $\tau$.

La démonstration de la proposition repose sur la correspondance de Riemann-Hilbert irrégulière, telle qu'elle est expliquée par exemple dans \cite{Deligne78} (voir aussi \cite{B-V89}, \cite[Chap\ptbl IV]{Malgrange91}). La structure formelle du supposé $(\cV_t,\nabla_t)$ étant fixée par la décomposition, il suffit, pour montrer son existence, de l'enrichir par une structure de Stokes (passer d'un système local $I$-gradué à un système local $I$-filtré, dans le langage de \loccit). Une fois la décomposition formelle fixée, la structure de Stokes dépend de la structure combinatoire sur le cercle unité des ensembles définis par les inéquations $\arg(t\ga-t\gb)\in[-\pi/2,\pi/2]$ pour $\ga\neq\gb$ intervenant dans la décomposition de $(\wh\cV,\wh\nabla)$. Ces intervalles étant indépendants de $t>0$, la structure de Stokes donnée par $(\cV,\nabla)$ pour $t=1$ se prolonge de manière unique pour tout $t>0$.

Pour adapter cette démonstration il convient, dans la situation locale considérée plus haut,
\begin{enumeratea}
\item\label{enum:hbvariablea}
d'étendre, pour un fibré à connexion méromorphe sur~$X$ à pôles le long de~$D$, qui est sauvage et bon, la théorie asymptotique de Sibuya-Majima; on a besoin pour cela de l'existence d'un bon réseau (voir \S\ref{subsec:goodwild}) pour argumenter par récurrence sur la dimension;

\item\label{enum:hbvariableb}
d'étendre à ce cas la correspondance de Riemann-Hilbert irrégulière; Mochizuki a ainsi redécouvert la notion de système local $I$-filtré de \cite{Deligne78}, l'a adaptée en dimension quelconque en prouvant aussi l'efficacité de cette approche, du point de vue métrique notamment;

\item\label{enum:hbvariablec}
d'appliquer, \emph{pour $\hb$ fixé}, l'analogue de la proposition \ref{prop:isomono} à $\ccP\cE^\hb$ avec le multiplicateur $1/(1+|\hb|^2)$ pour obtenir un prolongement noté $\ccQ\cE^\hb$;

\item\label{enum:hbvariabled}
d'étendre aussi la correspondance de Riemann-Hilbert irrégulière aux familles de connexions paramétrées par $\hb$, comme celles induites par une $\hb$-connexion si $\hb\neq0$;

\item\label{enum:hbvariablee}
de corriger la dépendance non holomorphe des facteurs exponentiels par un argument analogue à celui de la proposition \ref{prop:isomono} en choisissant convenablement le multiplicateur (correspondant à $t$) pour que, en fixant $\hb$ on retrouve $\ccQ\cE^\hb$; la caractérisation du prolongement $\ccQ\cE$ ainsi obtenu par une condition de croissance modérée nécessite de modifier la métrique $h$ et de la faire dépendre de $\hb$; [un tel type de correction apparaît déjà dans \cite{Szabo07} pour un calcul de transformé de Nahm d'un fibré de Higgs, et dans \cite{Bibi04} pour un calcul de transformé de Fourier];

\item\label{enum:hbvariablef}
de recoller la construction précédente, faite pour $\hb\neq0$, à $\ccP\cE^0$ (dans le cas Higgs, il n'y a pas de structure de Stokes et la décomposition \eqref{eq:decomp} est déjà holomorphe).
\end{enumeratea}
On en déduit:

\begin{theo}[{\cite[Th\ptbl11.12]{Mochizuki08}}]\label{th:QE}
Si $(E,\theta,h)$ est un fibré de Higgs harmonique sauvage et bon, il existe un unique $\cO_{X\times\CC_\hb}(*(D\times\CC_\hb))$-module localement libre~$\ccQ\cE$ à $\hb$\nobreakdash-connexion méromorphe dont la restriction à chaque $\hb$ soit égale à  $(\ccQ\cE^\hb,\nabla^\hb)$.
\end{theo}

\begin{rema}
Le fait que le phénomène de Stokes n'apparaisse pas pour certaines questions reliées au cas sauvage (il n'apparaît pas dans \cite{B-B03} par exemple) provient du fait que ce phénomène n'existe pas pour les fibrés de Higgs et que, pour le cas des fibrés méromorphes plats, plusieurs points peuvent se traiter à l'aide d'une approximation à un ordre assez grand de la structure formelle.
\end{rema}

\section{$\cD$-modules avec structure de twisteur sauvage}\label{sec:Dmodtw}

Nous avons vu à l'étape \eqref{enum:programme1} du \S\ref{sec:strategie} que le prolongement d'une variation de structure de twisteur polarisée sur $X^o$ va se chercher dans une catégorie de quadruplets $(\cM',\cM'',C,\cS)$, où $\cS$ sera la polarisation. La construction $\ccQ\cE$ du théorème \ref{th:QE} est l'étape principale pour construire $\cM'$, et on posera $\cM''=\cM'$ et $\cS=\id$ pour avoir un modèle simple quand le poids $w$ est nul. Pour passer de $\ccQ\cE$ à $\cM'$, il manque l'analogue du prolongement intermédiaire $j_{!*}$. Il aussi important de prolonger $C$ obtenu à partir de la métrique harmonique comme dans le dictionnaire du \S\ref{sec:dictionnaire}.

Un point essentiel à assurer, pour cette construction, et pour les constructions qui vont suivre, est que les $\cR_{X\times\CC_\hb}$-modules soient \emph{stricts}, c'est-à-dire sans $\cO_{\CC_\hb}$-torsion. On la sous-entendra cependant dans la suite. Par exemple, pour l'image directe par un morphisme propre, la préservation de cette condition est analogue à la dégénérescence en $E_1$ de la suite spectrale Hodge $\implique$ de~Rham.

\subsection{Les cycles proches}
Nous nous intéressons dans ce paragraphe à des propriétés locales pour un $\cR_{X\times\CC_\hb}$-module holonome.

La construction du foncteur des cycles proches modérés $\psi_f^\rmod$ par rapport à une fonction holomorphe $f:X\to\CC$ pour un $\cD_X$-module holonome $M$ s'appuie sur le théorème d'existence d'un polynôme de Bernstein-Sato, et s'exprime comme la graduation par rapport à la filtration de Kashiwara-Malgrange du $\cD_X$-module. Alors $\psi_f^\rmod M$ est un $\cD_X$-module holonome muni d'un endomorphisme semi-simple et d'un endomorphisme nilpotent $\rN$ qui commutent. Pour un $\cR_{X\times\CC_\hb}$-module holonome $\cM$, un tel polynôme de Bernstein-Sato n'existe pas nécessairement, mais on peut définir la sous-catégorie des $\cR_{X\times\CC_\hb}$-modules modérément spécialisables le long de tout germe de fonction holomorphe sur~$X$ en imposant l'existence d'une équation fonctionnelle de type Bernstein-Sato. Pour $\cM$ spécialisable, $\psi_f^\rmod\cM$ est défini, à support dans $\{f=0\}\times\CC_\hb$ et muni de deux endomorphismes comme plus haut.

Pour les $\cD_X$-modules holonomes à singularités irrégulières, ce foncteur $\psi_f^\rmod$ peut être de peu d'utilité. Par exemple, dans la situation de l'exemple \ref{exem:rankone}, pour toute fonction holomorphe $f(z)$ telle que $f(0)=0$, le foncteur $\psi_f^\rmod$ appliqué au $\cD_\Delta$-module $\cD_\Delta/\cD_\Delta\cdot\nobreak(z^2\partial_z+1)$ donne pour résultat~$0$. Dans \cite{Deligne83}, Deligne a défini à partir de~$\psi_f^\rmod$ un foncteur que nous notons $\psi_f^\Del$, et qui permet d'éviter ce comportement trivial: on a $\psi_f^\Del(\cM,\nabla)=\bigoplus_\ga\psi_f^\rmod(\cM,\nabla+\rd\ga)$, où $\ga$ parcourt l'ensemble des parties polaires de la variable $z^{1/q}$ et $q$ est un entier quelconque $\geq1$. On peut alors définir la sous-catégorie des $\cR_{X\times\CC_\hb}$-modules Deligne-spécialisables le long de tout germe de fonction holomorphe sur~$X$, et pour $\cM$ Deligne-spécialisable, $\psi_f^\Del\cM$ est à support dans $\{f=0\}\times\CC_\hb$ et est muni de deux endomorphismes comme plus haut. Dans ce cadre, la filtration \og monodromique\fg $\rM_\bbullet\psi_f^\Del\cM$ associée à l'endomorphisme nilpotent $\rN$ est bien définie (\cite[Prop\ptbl1.6.1]{Deligne80}), et on considérera ci-dessous les foncteurs $\gr_\ell^\rM\psi_f^\Del$. Ces foncteurs s'étendent de manière naturelle aux accouplements~$C$.

Enfin, la propriété de décomposabilité suivant le support, introduite par Saito, sera aussi importante. Les $\cR_{X\times\CC_\hb}$-modules holonomes stricts $\cM$ que nous considérons ont pour support un sous-ensemble analytique (ou algébrique) fermé $Z\times\CC_\hb$ de $X\times\CC_\hb$. Au voisinage de tout point~$x$ de~$Z$, on dispose des sous-modules $\cM_i$ à support dans un germe de sous-ensemble analytique fermé irréductible~$Z_i$ de~$Z$ en~$x$. La condition de \emph{S-décomposabilité} consiste à demander une décomposabilité locale $\cM=\bigoplus_i\cM_i$ en tout point de~$Z$. Si on dispose d'un accouplement~$C$, on impose aussi la diagonalité de~$C$ vis-à-vis de cette décomposition.

\subsection{La catégorie des $\cD$-modules holonomes avec structure de twisteur polarisable sauvage}

Cette catégorie est définie (voir \cite{Bibi06b}) en suivant le procédé de Saito pour les modules de Hodge polarisables, par récurrence sur la dimension du support.

\begin{defi}\label{def:wildtw}
La catégorie $\MTw_{\leq d}(X,w)$ des $\cD$-modules holonomes avec structure de twisteur sauvage pure de poids $w\in\ZZ$ est la sous-catégorie pleine de la catégorie des triplets $\cT=(\cM',\cM'',C)$ dont les objets satisfont aux propriétés suivantes:

\smallskip\noindent
\textup{(HSD)} $\cT$ est holonome, S-décomposable et a un support de dimension $\leq d$ dans~$X$.

\smallskip\noindent
$(\MTw_{>0})$ Pour tout ouvert $U\subset X$ et toute fonction holomorphe $f:U\to\nobreak\CC$, $\cT$ est Deligne-spécialisable le long de $\{f=0\}$ et, pour tout entier $\ell\geq0$, le triplet $\gr_\ell^{\rM}\Psi_f^\Del \cT$ est un objet de $\MTw_{\leq d-1}(X,w+\ell)$.

\smallskip\noindent
$(\MT_0)$
pour tout $x_o\in X$, la S-composante $(\cM'_{\{x_o\}},\cM''_{\{x_o\}},C_{\{x_o\}})$ est une \og masse de Dirac\fg en $x_o$ portant une structure de twisteur pure de poids $w$.
\end{defi}

On peut aussi définir les objets polarisables par des contraintes analogues sur la polarisation $\cS$.

Nous avons indiqué au \S\ref{sec:strategie}\eqref{enum:programme4} que le théorème \ref{th:HLT} s'applique à ces objets polarisables. Nous renvoyons à \cite[Chap\ptbl18]{Mochizuki08} pour la démonstration, qui s'inspire de la démonstration de Saito pour les modules de Hodge polarisables.

\subsection{Fin de la démonstration du théorème de Lefschetz difficile}

Soit $Z^o$ une variété quasi-projective lisse irréductible. Nous avons vu (\S\ref{sec:dictionnaire}) qu'une variation de structure de twisteur polarisable pure de poids $0$ sur $Z^o$ correspond à un fibré de Higgs harmonique. Nous dirons que cette variation est \emph{sauvage} s'il existe une compactification projective $Z'$ de $Z^o$ telle que $Z'\moins Z^o$ soit un diviseur à croisements normaux dont toutes les composantes sont lisses et que le fibré de Higgs harmonique soit sauvage et bon le long de ce diviseur. On note $\VTPw(Z^o,0)$ la catégorie correspondante.

Soit d'autre part $Z$ une compactification projective quelconque de $Z^o$ contenue dans une variété projective lisse $X$ et $\MTPw(Z,Z^o,0)$ la catégorie des $\cD_X$-modules holonomes avec structure de twisteur sauvage pure de poids $0$ et polarisable, lisses sur~$Z^o$ et sans composante à support dans un fermé strict de $Z$.

\begin{theo}[{\cite[Th\ptbl19.2]{Mochizuki07}, \cite[Cor\ptbl19.1.4]{Mochizuki08}}]\label{th:equiv}
La restriction à $Z^o$ définit un foncteur $\MTPw(Z,Z^o,0)\to\VTPw(Z^o,0)$. Ce foncteur est une équivalence de catégories.
\end{theo}

\begin{rema}[Comment changer de poids]
Les objets ci-dessus sont munis d'un twist de Tate $(k)$ pour tout $k\in\frac12\ZZ$, de manière analogue aux structures de Hodge complexes. Ce twist fait passer du poids $0$ au poids $-2k$.
\end{rema}

Ce théorème répond au point \eqref{enum:programme6} du \S\ref{sec:strategie}. Un des points difficiles est l'essentielle surjectivité du foncteur. La construction du $\cR_{X\times\bS}$-module $\ccQ\cE$, aboutissement du \S\ref{sec:prolong}, en est l'ingrédient essentiel. La construction du prolongement de l'accouplement $C$, ainsi que la démonstration des propriétés $\MTPw$ de l'objet $\cT$ ainsi obtenu, dans le cas où~$Z$ est lisse, $Z\moins Z^o$ est un diviseur à croisements normaux et la variation est sauvage et bonne le long de $Z\moins Z^o$, sont expliquées dans \hbox{\cite[Chap\ptbl12]{Mochizuki08}} (et~dans \cite[Chap\ptbl18]{Mochizuki07} pour le cas modéré). Enfin, le cas général est traité dans le chapitre~19 de \loccit, en utilisant, comme Saito le faisait dans \cite{MSaito86}, une version relative du théorème \ref{th:HLT} pour passer d'une bonne compactification de $Z^o$ à une moins bonne.

\section{Théorie de Hodge sauvage}

\subsubsection*{Variations de structure de Hodge complexe et variations intégrables de structure de twisteur}
La théorie précédente ne fournit que peu d'invariants numériques nouveaux pour les $\cD$-modules holonomes simples, contrairement à la théorie de Hodge classique qui fournit les nombres de Hodge $h^{p,q}$. On sait (\cite{Simpson97}) que les variations de structure de twisteur polarisable de poids $w$ qui proviennent d'une variation de structure de Hodge polarisable de poids $w$ sont celles qui sont munies d'une action naturelle de~$\CC^*$ (sur le facteur $\CC_\hb$), action qui permet de récupérer la graduation de Hodge. On peut aussi les caractériser comme celles admettant une action infinitésimale de $\CC^*$, lorsque la base $X^o$ est quasi-projective, si on impose la condition de modération à l'infini (théorème~\ref{th:HS} ci-dessous).

\begin{defi}
Une variation de structure de twisteur $(\cH',\cH'',C)$ (voir le \S\ref{sec:dictionnaire}) est \emph{intégrable} si les $\hb$-connexions sur $\cH'$ et $\cH''$ proviennent d'une connexion plate (absolue) avec un pôle de rang de Poincaré égal à $1$ le long de $\hb=0$, et si $C$ est compatible (en un sens naturel) à ces connexions.
\end{defi}

Autement dit, il existe $\nabla:\cH'\to\frac1\hb\Omega^1_{X\times\CC_\hb}(\log\{\hb\!=\!0\})\otimes\cH'$ (idem pour~$\cH''$), telle que $\nabla^2=0$, que la composante sur $\frac1\hb\Omega^1_{X\times\CC_\hb/\CC_\hb}$ soit la $\hb$-connexion, et enfin
\[
\hb\frac{\partial}{\partial\hb}\,C(m',\sigma^*\ov{m''})=C(\hb\nabla_{\partial_\hb}m',\sigma^*\ov{m''})-C(m',\sigma^*\ov{\hb\nabla_{\partial_\hb}m{}''}).
\]

\begin{theo}[voir {\cite[Th\ptbl6.2]{H-S08}}]\label{th:HS}
Les variations intégrables de structure de twisteur polarisée de poids $w$ sur un ouvert de Zariski $X^o$ de~$X$ (projective lisse) qui sont modérées à l'infini correspondent bijectivement aux variations de structure de Hodge complexe polarisée munies d'un automorphisme semi-simple auto-adjoint pour~$h$.
\end{theo}

Plus généralement, on peut alors considérer les variations intégrables et sauvages de structure de twisteur pure polarisée sur une variété quasi-projective comme des analogues irréguliers des variations de structure de Hodge complexe polarisée.

Hertling \cite{Hertling01} a observé qu'un nouvel invariant apparaît dans ce cadre, déjà considéré par les physiciens Cecotti et Vafa \cite{C-V91, C-F-I-V92}, appelé \og nouvel indice super-symétrique\fg. Une structure de twisteur pure de poids~$0$ et polarisée correspond simplement à un espace vectoriel complexe muni d'une forme hermitienne définie positive. Elle est intégrable si elle se trouve de plus munie de deux endomorphismes~$\cU$ et $\cQ$, avec~$\cQ$ auto-adjoint relativement à la forme hermitienne. L'espace vectoriel se décompose suivant les valeurs propres $p\in\RR$ de $\cQ$, et les composantes ont une dimension notée $h^{p,-p}$. Pour une structure de Hodge complexe polarisée pure de poids~$0$, ceci n'est autre que la décomposition de Hodge et les nombres de Hodge, avec $p\in\ZZ$ dans ce cas, et on a aussi $\cU=0$.

Dans une variation paramétrée par $x\in X$, l'exposant $p$ peut varier avec $x$ de manière analytique réelle, alors qu'il est constant (entier) pour les variations de structure de Hodge. Aussi, regrouper les espaces propres de $\cQ$ suivant les $p$ ayant même partie entière, pour avoir une décomposition de Hodge au sens usuel, peut provoquer des sauts de dimension suivant les valeurs de~$x$.

Le comportement de cet indice à l'infini d'une telle variation (modérée ou sauvage) est analysé dans \cite{Bibi08} en dimension $1$ et dans \cite{Mochizuki08b} en toute dimension.

\subsubsection*{Structures réelles et rationnelles}
Hertling \cite{Hertling01} a aussi considéré de telles variations avec une structure réelle (structure qu'il appelle TERP). Le théorème précédent est en fait montré avec structure réelle.

Plus récemment, Katzarkov, Kontsevich et Pantev \cite{K-K-P08} ont proposé la notion de structure rationnelle (voir aussi \cite{Bibi11}), et dans ce cadre l'appellation de \og variation de structure de Hodge non commutative\fg.

\subsubsection*{Théorie de Hodge mixte sauvage}
De manière analogue à la théorie des modules de Hodge mixtes de M\ptbl Saito \cite{MSaito87}, T\ptbl Mochizuki \cite{Mochizuki11} a développé la théorie des $\cD$-modules holonomes avec structure de twisteur mixte sauvage, éventuellement intégrable. Il explicite en particulier un foncteur de dualité, qui n'était pas défini dans le cadre précédent, qui permet notamment de définir la notion de structure réelle. Seule manque encore la structure rationnelle, mais les arguments de \cite{Mochizuki10} devraient pouvoir s'appliquer aussi à ce cadre.

\subsubsection*{Application à la cohomologie quantique}
Les résultats de H\ptbl Iritani sur la symétrie miroir pour les variétés de Fano toriques \cite{Iritani09,Iritani09b} (voir aussi \cite{R-S10}), joints aux résultats sur la transformation de Fourier de \cite{Bibi05}, permettent de montrer que le $\cD$\nobreakdash-module quantique d'une variété de Fano torique sous-tend une variation de structure de Hodge non commutative.

\backmatter
\newcommand{\SortNoop}[1]{}\def\cprime{$'$}
\providecommand{\bysame}{\leavevmode ---\ }
\providecommand{\og}{``}
\providecommand{\fg}{''}
\providecommand{\smfandname}{\&}
\providecommand{\smfedsname}{\'eds.}
\providecommand{\smfedname}{\'ed.}
\providecommand{\smfmastersthesisname}{M\'emoire}
\providecommand{\smfphdthesisname}{Th\`ese}

\end{document}